\pdfoutput=1

\documentclass[10pt,a4paper]{amsart}

	\usepackage[utf8]{inputenc} 
	\usepackage[T1]{fontenc} 
	\usepackage[british]{babel} 
	
	\usepackage[backend=bibtex,style=alphabetic,defernumbers=true,firstinits=true,doi=false,isbn=false,url=false,useprefix=true,maxbibnames=9,maxcitenames=9,sorting = nyt]{biblatex} 
	\bibliography{ARTICLE.bib} 
	\DeclareFieldFormat{labelalpha}{\mkbibbold{#1}}
	\DeclareFieldFormat{extraalpha}{\mkbibbold{\mknumalph{#1}}}
	\AtBeginBibliography{
		\DeclareFieldFormat{labelalpha}{#1}
		\DeclareFieldFormat{extraalpha}{\mknumalph{#1}}}
	
	\usepackage{geometry}
	\usepackage{amsmath}
	\usepackage{amssymb}
	\usepackage{amsthm}
	\usepackage{mathtools}  
	\usepackage{tikz-cd} 
	\usepackage[mathcal]{euscript} 
	\usepackage{mathrsfs} 
	
	\usepackage{enumitem} 
	
	\definecolor{RadboudRed}{RGB}{190, 49, 26} 
	\definecolor{RadboudRedDark}{RGB}{144, 36, 22}
	\usepackage[colorlinks=true,linkcolor=RadboudRed,citecolor={green!50!black},urlcolor=orange]{hyperref} 
	\usepackage[capitalise,nameinlink]{cleveref}
	\crefdefaultlabelformat{#2\textbf{#1}#3} 

	\newcommand{\demph}[1]{{\emph{#1}}}    
	
	\newcommand{\Cstar}{C$^\ast$}	
	\let\im\relax\DeclareMathOperator{\im}{im}	
	\newcommand{\id}{\mathrm{id}}			
	\newcommand{\pr}{\mathrm{pr}}			
	\DeclareMathOperator{\Hom}{Hom}			
	
	\newcommand{\ar}{\mathrm{ar}}	
	\DeclareMathOperator{\dom}{dom}
	
	\newcommand{\catfont}[1]{\mathbf{#1}}

	\newcommand{\Mnfd}{\catfont{Mnfd}}
	
	\newcommand{\Diff}{\catfont{Diffeol}}

	\newcommand{\Germ}{\catfont{Germ}}		

	\newcommand{\Act}{\catfont{Act}}
	\newcommand{\DiffBiBund}{\catfont{DiffeolBiBund}}

	\newcommand{\source}{\mathrm{src}} 
	\newcommand{\target}{\mathrm{trg}} 
	\newcommand{\inverse}{\mathrm{inv}} 
	
	\newcommand{\unit}{\mathrm{u}} 
	\newcommand{\grpd}[1]{{#1}\rightrightarrows{#1}_0}
	
	\newcommand{\laction}{\curvearrowright}	
	\newcommand{\raction}{\curvearrowleft}
	\newcommand{\LAalong}[1]{{\laction\hspace{-5pt}^{#1}\hspace{1pt}}}
	\newcommand{\RAalong}[1]{{\hspace{-1pt}~^{#1}\hspace{-5pt}\raction}}
	\newcommand{\Orb}{\mathrm{Orb}}   
	
	\newcommand{\diffeology}{\mathcal{D}}
	\newcommand{\Param}{\mathrm{Param}}
	
	\renewcommand{\geq}{\geqslant}
	
	\newcommand{\rra}{\rightrightarrows}
	
	\newcommand{\mc}[1]{\mathcal{#1}}
	\newcommand{\const}{\mathrm{const}}
	\newcommand{\To}{\Rightarrow}
	\newcommand{\tensor}{\otimes}
	
	\newcommand{\diff}{\mathrm{Diff}}
	\newcommand{\comp}{\mathrm{comp}}

	\theoremstyle{plain}
	\newtheorem{theorem}{Theorem}[section]
	\newtheorem*{theorem*}{Theorem}
	\newtheorem{proposition}[theorem]{Proposition}
	\newtheorem{corollary}[theorem]{Corollary}
	\newtheorem*{conjecture*}{Conjecture}

	\newtheorem{lemma}[theorem]{Lemma}
	\newtheorem*{lemma*}{Lemma}
	\newtheorem{question}[theorem]{Question}
	\newtheorem*{main_theorem_introduction}{\cref{theorem:weakly invertible bibundles are the biprincipal ones}} 

	\theoremstyle{definition}
	\newtheorem{definition}[theorem]{Definition}
	\newtheorem{example}[theorem]{Example}

	\newtheorem{construction}[theorem]{Construction}
	
	\crefname{wishlist}{Wishlist}{Wishlists}
	\crefname{part}{Chapter}{Chapters}
	\crefname{construction}{Construction}{Constructions}

\title{\textbf{Diffeological Morita Equivalence}\vspace{-\baselineskip}}
\author{Nesta van der Schaaf}
\thanks{\emph{E-mail address:} \url{n.schaaf@ed.ac.uk} or \url{nestavanderschaaf@gmail.com}}
\date{\today}

\begin{document}
		\begin{abstract}
			We introduce a new notion of \emph{Morita equivalence} for \emph{diffeological groupoids}, generalising the original notion for Lie groupoids. For this we develop a theory of \emph{diffeological groupoid actions}, \emph{-bundles} and \emph{-bibundles}. We define a notion of \emph{principality} for these bundles, which uses the notion of a \emph{subduction}, generalising the notion of a Lie group(oid) principal bundle. We say two diffeological groupoids are Morita equivalent if and only if there exists a \emph{biprincipal} bibundle between them. Using a Hilsum-Skandalis tensor product, we further define a composition of diffeological bibundles, and obtain a bicategory $\DiffBiBund$. Our main result is the following: a bibundle is biprincipal if and only if it is \emph{weakly invertible} in this bicategory. This generalises a well known theorem from the Lie groupoid theory. As an application of the framework, we prove that the \emph{orbit spaces} of two Morita equivalent diffeological groupoids are diffeomorphic. We also show that the property of a diffeological groupoid to be \emph{fibrating}, and its \emph{category of actions}, are Morita invariants.%
			
		\end{abstract}
		\maketitle
\section{Introduction}
\label{section:introduction}
\emph{Diffeology} originates from the work of J.-M. Souriau \cite{souriau1980groupes,souriau1984groupes} and his students \cite{donato1983exemple,donato1984revetements,iglesias1985fibrations} in the 1980s. The main objects of this theory are \emph{diffeological spaces}, a type of generalised smooth space that extends the traditional notion of a smooth manifold. They make for a convenient framework that deals well with (singular) quotients, function spaces (or otherwise infinite-dimensional objects), fibred products (or otherwise singular subspaces), and other constructions that lie beyond the realm of classical differential topology. As many of these constructions naturally occur in differential topology and -geometry, and since they cannot be studied with their standard tools, diffeology has become a useful addition to the geometer's toolbox. 

\emph{Diffeological groupoids} have recently garnered attention in the mathematical physics of general relativity \cite{blohmann2013groupoid,glowacki2019}, foliation theory \cite{androulidakis2019diffeological,garmendia2019hausdorff,macdonald2020holonomygroupoids}, the theory of algebroids \cite{androulidakis2020integration}, the theory of (differentiable) stacks \cite{roberts2018smooth,watts2019diffeological}, and even in relation to noncommutative geometry \cite{iglesias2018noncommutative,iglesias2020quasifolds}. In all but one of these fields (general relativity), the notion of \emph{Morita equivalence} is an important one. Yet, as the authors of \cite[p.3]{garmendia2019hausdorff} point out: ``The theory of Morita equivalence for diffeological groupoids
has not been developed yet.'' In the current paper we present one possible development of such a notion, based on the results of the author's Master thesis \cite{schaaf2020diffeology-groupoids-and-ME}. This development is a generalisation of the theory of Hilsum-Skandalis bibundles and the Morita equivalence of Lie groupoids, where many definitions and proofs, and certainly the general idea, extend quite straightforwardly to the diffeological case. The main exception is that we need to replace surjective submersions with so-called \emph{subductions}. This special type of smooth map is, even on smooth manifolds, slightly weaker than the notion of a surjective submersion, but it turns out that they still share enough of their properties so that the entire theory can be developed%
\footnote{This is essentially due to the fact that the subductions are the \emph{strong epimorphisms} in the category of diffeological spaces \cite[Proposition 5.10]{baez2011convenient}.}. 
This development proceeds roughly as follows: based on the notions of \emph{actions} and \emph{bundles} defined in \cref{section:diffeological groupoid actions and bundles}, we define a diffeological version of a \emph{bibundle} between groupoids (\cref{definition:diffeological groupoid bibundles}). These stand in analogy to \emph{bimodules} for rings, and can be treated as a generalised type of morphism between groupoids. This gives a \emph{bicategory} $\DiffBiBund$ of diffeological groupoids, bibundles, and \emph{biequivariant maps} (\cref{theorem:bicategory DiffBiBund}).
Using the aforementioned notion of a \emph{subduction} (\cref{definition:subduction}), we define \emph{biprincipality} of bibundles, and with this, we obtain a notion of \emph{Morita equivalence} for diffeological groupoids (\cref{definition:Morita equivalence and biprincipality}). In the bicategory we also get a notion of equivalence, by way of the \emph{weak isomorphisms}. A morphism in a bicategory is called \emph{weakly invertible} if it is invertible \emph{up to 2-isomorphism}. Two objects in a bicategory are called \emph{weakly isomorphic} if there exists a weakly invertible morphism between them. The main point of this paper is to prove a \emph{Morita theorem} for diffeological groupoids, characterising the weakly invertible bibundles, and hence realising Morita equivalence as a particular instance of weak isomorphism:
	\begin{main_theorem_introduction}[Morita theorem]
		A diffeological bibundle is weakly invertible if and only if it is biprincipal. In other words, two diffeological groupoids are Morita equivalent if and only if they are weakly isomorphic in the bicategory $\DiffBiBund$. 
	\end{main_theorem_introduction}
A Morita theorem for Lie groupoids has been known in the literature for some time, see e.g. \cite[Proposition 4.21]{landsman2001quantized}. Throughout the paper, we shall point out some differences between the diffeological- and Lie theories. The main difference is that, due to technical constraints, a Morita theorem for Lie groupoids only holds in the restricted setting of \emph{left principal} bibundles. The main improvement of \cref{theorem:weakly invertible bibundles are the biprincipal ones} over the classical Lie Morita theorem, besides the generalisation to diffeology, is therefore that it considers also a more general class of bibundles. Besides this improvement, with this paper we hope to contribute a complete account of the basic theory of bibundles and Morita equivalence of groupoids, providing detailed proofs and constructions of most necessary technical results, and culminating in a proof of the main \cref{theorem:weakly invertible bibundles are the biprincipal ones}. A brief outline of the contents of the paper is as follows.

We briefly recall the definition of a diffeology in \cref{section:diffeology}. In particular, we describe the diffeologies of fibred products (pullbacks) and quotients, since they will be important to describe the smooth structure of the orbit space and space of composable arrows of a groupoid. We also define and study the behaviour of \emph{subductions}, especially in relation to fibred products.

In \cref{section:diffeological groupoids} we define \emph{diffeological groupoids}, and highlight some examples from the literature.

\cref{section:diffeological groupoid actions and bundles,section:diffeological bibundles} contain the main contents of this paper. In them, we define the notions of smooth groupoid \emph{actions} and \emph{-bundles}. For the latter we give a new notion of \emph{principality}, generalising the notion of a principal Lie group(oid) bundle. This leads naturally to the definition of a \emph{biprincipal bibundle}, and hence to our definition of \emph{Morita equivalence}. The remainder of \cref{section:diffeological bibundles} is dedicated to a proof of \cref{theorem:weakly invertible bibundles are the biprincipal ones}.

In \cref{section:some applications}, we describe some \emph{Morita invariants}, by generalising some well known theorems from the Lie theory. We prove: the property of a diffeological groupoid to be \emph{fibrating} is preserved under our notion of Morita equivalence; the \emph{orbit spaces} of two Morita equivalent diffeological groupoids are diffeomorphic; and the categories of representations of two Morita equivalent diffeological groupoids are categorically equivalent. 

Lastly, in \cref{section:closing section}, we discuss the question of diffeological Morita equivalence between Lie groupoids. We end the paper with the open \cref{question:does inclusion pseudofunctor reflect weak equivalence}, and some suggestions for future research.

\smallskip
\noindent \textbf{Acknowledgements.} The author thanks Klaas Landsman and Ioan M\u{a}rcu\textcommabelow{t} for being the supervisor and second reader of his Master thesis, respectively, and for encouraging him to write the current paper. He also thanks Klaas for feedback on the paper, and Patrick Iglesias-Zemmour for email correspondence.

\section{Diffeology}
	\label{section:diffeology}
	One of the main conveniences of \emph{diffeology}%
	\footnote{The etymology of the word is explained in the afterword to \cite{iglesias2013diffeology}. Souriau first used the term \emph{``diff\'erentiel''}, as in `differential' (from the Latin \emph{differentia}, ``difference''). Through a suggestion by Van Est, the name was later changed to \emph{``diff\'eologie,''} as in \emph{``topologie''} (`topology', from the Ancient Greek \emph{t\'opos}, ``place,'' and \emph{-(o)logy}, ``study of''). Hence the term: diffeology.}
	is that the category $\Diff$ of diffeological spaces and smooth maps (\cref{definition:diffeology}) is complete, cocomplete, (locally) Cartesian closed, and in fact a quasitopos \cite[Theorem 3.2]{baez2011convenient}. This means that we can perform many categorical constructions that are unavailable in the category $\Mnfd$ of smooth manifolds. From these, the ones that are important for us are pullbacks and quotients. We discuss both of these explicitly below. The approach of diffeology has been compared to other theories of generalised smooth spaces in \cite{stacey2011comparative,batubenge2017diffeologicalfrolicher}. For some historical remarks we refer to \cite{iglesias2013beginning,iglesias2019introduction} and \cite[Chapter I]{schaaf2020diffeology-groupoids-and-ME}. The main reference for this section is the textbook \cite{iglesias2013diffeology} by Iglesias-Zemmour, in which nearly all of the theory below is already developed.

	\begin{definition}
		A \emph{Euclidean domain} is an open subset $U\subseteq \mathbb{R}^m$, for arbitrary $m\in\mathbb{N}_{\geq 0}$. A \demph{parametrisation} on an arbitrary set $X$ is a function $U\to X$ defined on a Euclidean domain. We denote by $\Param(X)$ the set of all parametrisations on $X$.
	\end{definition}
	
	The basic idea behind diffeology is that it determines which parametrisations are \emph{`smooth'}, in such a way that it captures the properties of ordinary smooth functions on smooth manifolds. The precise definition is as follows:
	
	\begin{definition}[Axioms of Diffeology]
		\label{definition:diffeology}
		Let $X$ be a set. A \demph{diffeology} on $X$ is a collection of parametrisations $\diffeology_X\subseteq \Param(X)$, containing what we call \demph{plots}, satisfying the following three axioms:
		\begin{itemize}
			\item \emph{(Covering)} Every constant parametrisation $U\to X$ is a plot.
			\item \emph{(Smooth Compatibility)} For every plot $\alpha:U_\alpha\to X$ in $\diffeology_X$ and every smooth function $h:V\to U_\alpha$ between Euclidean domains, we have that $\alpha\circ h \in\diffeology_X$.
			\item \emph{(Locality)} If $\alpha:U_\alpha\to X$ is a parametrisation, and $(U_i)_{i\in I}$ an open cover of $U_\alpha$ such that each restriction $\alpha|_{U_i}$ is a plot of $X$, then $\alpha\in\diffeology_X$. 
		\end{itemize}
		A set $X$, paired with a diffeology: $(X,\diffeology_X)$, is called a \demph{diffeological space}. Although, usually we shall just write $X$.
		
		A function $f:(X,\diffeology_X)\to (Y,\diffeology_Y)$ between diffeological spaces is called \emph{smooth} if for every plot $\alpha\in\diffeology_X$ of $X$, the composition $f\circ\alpha\in\diffeology_Y$ is a plot of $Y$. The set of all smooth functions between such diffeological spaces is denoted $C^\infty(X,Y)$, and smoothness is preserved by composition. The category of diffeological spaces and smooth maps is denoted by $\Diff$, and the isomorphisms in this category are called \emph{diffeomorphisms}.
	\end{definition}
	
	\begin{example}
		\label{example:Euclidean and manifold diffeologies}
		Any Euclidean domain $U$ gets a canonical diffeology $\diffeology_U$, called the \emph{Euclidean diffeology}. Its plots are the parametrisations that are smooth in the ordinary sense of the word. Similarly, we get a canonical diffeology $\diffeology_M$ for any smooth manifold $M$, called the \emph{manifold diffeology}. With respect to these diffeologies, the notion of smoothness defined in \cref{definition:diffeology} agrees with the ordinary one. Hence the inclusion functor $\Mnfd\hookrightarrow\Diff$ is fully faithful, and we can adopt the previous definition without causing any confusion.
	\end{example}
	
	\begin{example}
		\label{example:coarse and discrete diffeologies}
		Any set $X$ carries two canonical diffeologies. First, the largest diffeology, $\diffeology^\bullet_X:=\Param(X)$, called the \emph{coarse diffeology}, containing all possible parametrisations. Letting $X^\bullet$ denote the diffeological space with the coarse diffeology, it is easy to see that every function $Z\to X^\bullet$ is smooth.
		On the other hand, the smallest diffeology on $X$ is $\diffeology_X^\circ$, containing all locally constant parametrisations. This is called the \emph{discrete diffeology}. Similar to the above, we find that every function $X^\circ\to Y$ is smooth.
	\end{example}
	
	\begin{example}
		For any two diffeological spaces $X$ and $Y$, there is a natural diffeology on the space of smooth functions $C^\infty(X,Y)$ called the \emph{standard functional diffeology} \cite[Article 1.57]{iglesias2013diffeology}. It is the smallest diffeology that makes the evaluation map $(f,x)\mapsto f(x)$ smooth. With these diffeologies, $\Diff$ becomes Cartesian closed.
	\end{example}
	
	\subsection{Generating families}
	\label{section:generating families}
	The Axiom of Locality in \cref{definition:diffeology} ensures that the smoothness of a parametrisation, or of a function between diffeological spaces, can be checked locally. This allows us to introduce the following notions, which will help us study interesting constructions, and will often simplify proofs.
	
	\begin{definition}
		\label{definition:generating families}
		Consider a family $\mc{F}\subseteq \Param(X)$ of parametrisations on $X$. There exists a smallest diffeology on $X$ that contains $\mc{F}$. We denote this diffeology by $\langle \mc{F}\rangle$, and call it the \emph{diffeology generated by $\mc{F}$}. If $\diffeology_X=\langle \mc{F}\rangle$, we say $\mc{F}$ is a \emph{generating family} for $\diffeology_X$. The elements of $\mc{F}$ are called \emph{generating plots}.
	\end{definition}
	
	The plots of the diffeology generated by $\mc{F}$ are characterised as follows: a parametrisation $\alpha:U_\alpha\to X$ is a plot in $\langle\mc{F}\rangle$ if and only if $\alpha$ is locally either constant, or factors through elements of $\mc{F}$. Concretely, this means that for all $t\in U_\alpha$ there exists an open neighbourhood $t\in V\subseteq U_\alpha$ such that $\alpha|_V$ is either constant, or of the form $\alpha|_V=F\circ h$, where $F:W\to X$ is an element in $\mc{F}$, and $h:V\to W$ is a smooth function between Euclidean domains. When the family $\mc{F}$ is \emph{covering}, in the sense that $\bigcup_{F\in\mc{F}}\im(F)=X$, then the condition for $\alpha|_V$ to be constant becomes redundant, and the plots in $\langle\mc{F}\rangle$ are locally just of the form $\alpha|_V=F\circ h$. 
	
	The main use of this construction is that we may encounter families of parametrisations that are not quite diffeologies, but that contain functions that we nevertheless want to be smooth. On the other hand, calculations may sometimes be simplified by finding a suitable generating family for a given diffeology. This simplification lies in the following result, saying that smoothness has only to be checked on generating plots:
	\begin{proposition}
		\label{proposition:smoothness defined by generating plots}
		Let $f:X\to Y$ be a function between diffeological spaces, such that $\diffeology_X$ is generated by some family $\mc{F}$. Then $f$ is smooth if and only if for all $F\in\mc{F}$ we have $f\circ F\in\diffeology_Y$. 
	\end{proposition}
	
	\begin{example}
		\label{example:wire diffeology}
		The \emph{wire diffeology} (called the \emph{spaghetti diffeology} by Souriau) is the diffeology $\diffeology_\mathrm{wire}$ on $\mathbb{R}^2$ generated by $C^\infty(\mathbb{R},\mathbb{R}^2)$. The resulting diffeological space is not diffeomorphic to the ordinary $\mathbb{R}^2$, since the identity map $\id_{\mathbb{R}^2}:(\mathbb{R}^2,\diffeology_{\mathbb{R}^2})\to (\mathbb{R}^2,\diffeology_\mathrm{wire})$ is not smooth.
	\end{example}
	
	\begin{example}
		\label{example:manifold diffeology is generated by atlas}
		The charts of a smooth atlas on a manifold define a generating family for the manifold diffeology from \cref{example:Euclidean and manifold diffeologies}. Since a manifold may have many atlases, this shows that similarly any diffeology may have many generating families.
	\end{example}
	\subsection{Quotients} 
	\label{section:quotients}
	We use the terminology from \cref{section:generating families} to define a natural diffeology on a quotient $X/{\sim}$. This question relates to a more general one: given a function $f:X\to Y$, and a diffeology $\diffeology_X$ on the domain, what is the smallest diffeology on $Y$ such that $f$ remains smooth? The following provides an answer:
	
	\begin{definition}
		\label{definition:pushforward diffeology}
		Let $f:X\to Y$ be a function between sets, and let $\diffeology_X$ be a diffeology on $X$. The \emph{pushforward diffeology} on $Y$ is the diffeology $f_\ast(\diffeology_X):= \langle f\circ\diffeology_X\rangle$, where $f\circ \diffeology_X$ is the family of parametrisations of the form $f\circ\alpha$, for $\alpha\in\diffeology_X$. The pushforward diffeology is the smallest diffeology on $Y$ that makes $f$ smooth.
	\end{definition}
	
	We can now use this to define a natural diffeology on a quotient space:
	\begin{definition}
		\label{definition:quotient diffeology}
		Let $X$ be a diffeological space, and let $\sim$ be an equivalence relation on the set $X$. We denote the equivalence classes by $[x]:=\lbrace y\in X:x\sim y\rbrace$. The \emph{quotient} $X/{\sim}$ is the collection of all equivalence classes, and comes with a \emph{canonical projection map} $p:X\to X/{\sim}$, which sends $x\mapsto [x]$. The \emph{quotient diffeology} on $X/{\sim}$ is defined as the pushforward diffeology $p_\ast(\diffeology_X)$ of $\diffeology_X$ along the canonical projection map. Naturally, with respect to this diffeology, the canonical projection map becomes smooth. 
	\end{definition}
	
	The quotient diffeology will be used extensively, where the equivalence relation will often be defined by the orbits of a group(oid) action, or as the fibres of some smooth surjection. The existence of the quotient diffeology for arbitrary quotients should be contrasted to the situation for smooth manifolds, where quotients often carry no natural differentiable structure at all, but where instead one could appeal to the \emph{Godement criterion} (\cite[Theorem 2, p. 92]{serre1965lie}). The following is an example of a quotient that does not exist as a smooth manifold, but whose diffeological structure is still quite rich:
	
	\begin{example}
		\label{example:irrational torus}
		The \emph{irrational torus} is the diffeological space defined by the quotient of $\mathbb{R}$ by an additive subgroup: $T_\theta:= \mathbb{R}/(\mathbb{Z}+\theta\mathbb{Z})$,
		where $\theta\in\mathbb{R}\setminus\mathbb{Q}$ is an arbitrary irrational number. Equivalently, it can be described as the leaf space of the Kronecker foliation on the 2-torus with irrational slope. The topology of this quotient contains only the two trivial open sets, yet its quotient diffeology is non-trivial\footnote{This shows that there are meaningful notions of smooth space that do not rely on the regnant philosophy of ``smooth space = topological space + extra structure.''}. They were first classified in \cite{donato1983exemple}, whose result is (amazingly) directly analogous to the classification of the irrational rotation algebras \cite{rieffel1981cstar}. This example is treated in detail in \cite[Section 2.3]{schaaf2020diffeology-groupoids-and-ME}.
	\end{example}
	
	\subsection{Fibred products}
	\label{section:fibred products}
	The second construction we need is that of \emph{fibred products}, which are the pullbacks in the category $\Diff$. Recall that if $f:X\to Z$ and $g:Y\to Z$ are two functions between sets with a common codomain, then the fibred product of sets is (up to unique bijection)
	\begin{equation*}
		X\times_Z^{f,g}Y
		:=
		\lbrace
		(x,y)\in X\times Y
		:
		f(x)=g(y)
		\rbrace.
	\end{equation*}
	When each set is equipped with a diffeology, we shall construct a diffeology on the fibred product in two steps. First we describe a natural diffeology on the product $X\times Y$, and then show how this descends to a diffeology on the fibred product as a subset. 
	
	\begin{definition}
		\label{definition:product diffeology}
		Let $X$ and $Y$ be two diffeological spaces. The \emph{product diffeology} on the Cartesian product $X\times Y$ is defined as 
		\begin{equation*}
			\diffeology_{X\times Y}
			:=
			\langle
			\diffeology_X\times \diffeology_Y
			\rangle,
		\end{equation*}
		where $\diffeology_X\times\diffeology_Y$ is the family of parametrisations of the form $\alpha_1\times \alpha_2$, for $\alpha_1\in\diffeology_X$ and $\alpha_2\in\diffeology_Y$. The plots in $\diffeology_{X\times Y}$ are exactly the parametrisations $\alpha:U_\alpha\to X\times Y$ such that $\pr_1\circ \alpha$ and $\pr_2\circ\alpha$ are plots of $X$ and $Y$, respectively. We assume that products are always furnished with their product diffeologies.
	\end{definition}
	
	It is clear that both projection maps $\pr_1$ and $\pr_2$ are smooth with respect to the product diffeology. The smooth functions into $X\times Y$ behave exactly as one would expect, where $f:A\to X\times Y$ is smooth if and only if the components $f_1=\pr_1\circ f$ and $f_2=\pr_2\circ f$ are smooth. 
	
	Next we define how the diffeology on a set $X$ transfers to any of its subsets:
	
	\begin{definition}
		\label{definition:subset diffeology}
		Consider a diffeological space $X$, and an arbitrary subset $A\subseteq X$. Let $i_A:A\hookrightarrow X$ denote the natural inclusion map. The \emph{subset diffeology} on $A$ is defined as
		\begin{equation*}
			\diffeology_{A\subseteq X}
			:=
			\lbrace
			\alpha\in\Param(A)
			:
			i_A\circ \alpha\in\diffeology_X
			\rbrace.
		\end{equation*}
		That is, $\alpha$ is a plot of $A$ if and only if when seen as a parametrisation of $X$, it is also a plot. We assume that a subset of a diffeological space is always endowed with its subset diffeology. 
	\end{definition}
	
	Since the fibred product $X\times_Z^{f,g}Y$ is a subset of the product $X\times Y$, the following definition is a natural combination of \cref{definition:product diffeology,definition:subset diffeology}:
	\begin{definition}
		\label{definition:fibred product diffeology}
		Let $f:X\to Z$ and $g:Y\to Z$ be two smooth maps between diffeological spaces. The \emph{fibred product diffeology} $\diffeology_{X\times_Z^{f,g}Y}$ on the set $X\times_Z^{f,g}Y$ is the subset diffeology it gets from the product diffeology on $X\times Y$. Concretely:
		\begin{equation*}
			\diffeology_{X\times_Z^{f,g}Y}
			=
			\lbrace
			\alpha\in\diffeology_{X\times Y}
			:
			f\circ\alpha_1
			=
			g\circ\alpha_2
			\rbrace.
		\end{equation*}
		That is, the plots of the fibred product are just plots of $X\times Y$, whose components satisfy an extra condition. We assume that all fibred products are equipped with their fibred product diffeologies.
	\end{definition}

	\subsection{Subductions} 
	Subductions are a special class of smooth functions that generalise the notion of surjective submersion from the theory of smooth manifolds. Since there is no unambiguous notion of tangent space in diffeology (cf. \cite{christensen2016tangent}), the definition looks somewhat different. For (more) detailed proofs of the results in this section, we refer to \cite[Article 1.46]{iglesias2013diffeology} and surrounding text, and \cite[Section 2.6]{schaaf2020diffeology-groupoids-and-ME}.
	
	\begin{definition}
		\label{definition:subduction}
		A surjective function $f:X\to Y$ between diffeological spaces is called a \emph{subduction} if $f_\ast(\diffeology_X)=\diffeology_Y$. Note that subductions are automatically smooth. 
	\end{definition}
	
	In the case that $f$ is a subduction, since it is then particularly a surjection, the family of parametrisations $f\circ\diffeology_X$ is covering, and hence the plots of $\diffeology_Y$ are all locally of the form $f\circ\alpha$, where $\alpha\in\diffeology_X$. In other words, $f$ is a subduction if and only if $f$ is smooth and the plots of $Y$ can locally be \emph{lifted} along $f$ to plots of $X$:
	\begin{lemma}
		\label{lemma:characterisation of subductions}
		Let $f:X\to Y$ be a function between diffeological spaces. Then $f$ is a subduction if and only if the following two conditions are satisfied:
		\begin{enumerate}
			\item The function $f$ is smooth.
			
			\item For every plot $\alpha:U_\alpha\to Y$, and any point $t\in U_\alpha$, there exists an open neighbourhood $t\in V\subseteq U_\alpha$ and a plot $\beta:V\to X$, such that $\alpha|_V=f\circ \beta$. 
		\end{enumerate}
	\end{lemma}
	Since many of the functions we encounter will naturally be smooth already, the notion of subductiveness is effectively captured by condition \emph{(2)} in this lemma. This can also be seen in the following simple example:
	
	\begin{example}
		\label{example:projection maps are subductions}
		Consider the product $X\times Y$ of two diffeological spaces $X$ and $Y$. The projection maps $\pr_1$ and $\pr_2$ are both subductions.
	\end{example}
	
	\begin{example}
		\label{example:quotient by surjection or subduction}
		For a surjective function $\pi:X\to B$ we get an equivalence relation on $X$, where two points are identified if and only if they inhabit the same $\pi$-fibre. The equivalence classes are exactly the $\pi$-fibres themselves. We denote the quotient set of this equivalence relation by $X/{\pi}$, and equip it with the quotient diffeology whenever $X$ is a diffeological space. If $\pi$ is a subduction, then there is a diffeomorphism $B\cong X/{\pi}$ \cite[Article 1.52]{iglesias2013diffeology}.
	\end{example}

	For subsequent use, we state here some useful properties of subductions with respect to composition:
	
	\begin{lemma}
		\label{lemma:properties of subductions}
		We have the following properties for subductions:
		\begin{enumerate}
			\item If $f$ and $g$ are two subductions, then the composition $f\circ g$ is a subduction as well.
			
			\item Let $f:Y\to Z$ and $g:X\to Y$ be two smooth maps such that the composition $f\circ g$ is a subduction. Then so is $f$. 
			
			\item Let $\pi:X\to B$ be a subduction, and $f:B\to Y$ an arbitrary function. Then $f$ is smooth if and only if $f\circ\pi$ is smooth. In fact, $f$ is a subduction if and only if $f\circ\pi$ is a subduction.
		\end{enumerate}
		
		\begin{proof}
			\emph{(1)} This is \cite[Article 1.47]{iglesias2013diffeology}.
			
			\emph{(2)} Assume $f:Y\to Z$ and $g:X\to Y$ are smooth, such that $f\circ g$ is a subduction. Take a plot $\alpha:U_\alpha\to Z$. Since the composition is a subduction, for every $t\in U_\alpha$ we can find an open neighbourhood $t\in V\subseteq U_\alpha$ and a plot $\beta:V\to X$ such that $\alpha|_V=(f\circ g)\circ \beta$. Since $g$ is smooth, we get a plot $g\circ\beta\in\diffeology_Y$, which is a local lift of $\alpha$ along $f$. The result follows by \cref{lemma:characterisation of subductions}.
			
			\emph{(3)} If $f$ is smooth, it follows immediately that $f\circ \pi$ is smooth. Suppose now that $f\circ \pi$ is smooth. We need to show that $f$ is smooth. For that, take a plot $\alpha:U_\alpha\to X$. Since $\pi$ is a subduction, we can find an open cover $(V_t)_{t\in U_\alpha}$ of $U_\alpha$ together with a family of plots $\beta_t:V_t\to X$ such that $\alpha|_{V_t}=\pi\circ\beta_t$. It follows that each restriction $f\circ\alpha|_{V_t}=f\circ\pi\circ\beta_t$ is smooth, and by the Axiom of Locality it follows that $f\circ\alpha\in\diffeology_Y$, and hence that $f$ is smooth. The claim about when $f$ is a subduction follows from \emph{(2)}.
		\end{proof}
	\end{lemma}
	
	We also collect the following noteworthy claim:
	\begin{proposition}[{\cite[Article 1.49]{iglesias2013diffeology}}]
		\label{proposition:injective subduction is diffeomorphism}
		An injective subduction is a diffeomorphism.
	\end{proposition}
	
	We recall now some elementary results on the interaction between subductions and fibred products, as obtained in \cite[Section 2.6]{schaaf2020diffeology-groupoids-and-ME}. We point out that if $f$ is a subduction, an arbitrary restriction $f|_A$ may no longer be a subduction. We know from \cref{example:projection maps are subductions} that the second projection map $\pr_2$ of a product $X\times Y$ is a subduction, but it is not always the case that the restriction of this projection to a fibred product $X\times_Z^{f,g}Y$ is a subduction as well. The following result shows that, to ensure this, it suffices to assume that $f$ is a subduction:
	\begin{lemma}
		\label{lemma:restriction of projection is subduction}
		Let $f:X\to Z$ be a subduction, and let $g:Y\to Z$ be a smooth map. Then the restricted projection map
		\begin{equation*}
			\left.\pr_2\right|_{X\times_Z^{f,g}Y}
			:
			X\times_Z^{f,g}Y
			\longrightarrow
			Y
		\end{equation*}
		is also a subduction. In other words, in $\Diff$, subductions are preserved under pullback.
		
		\begin{proof}
			Consider a plot $\alpha:U_\alpha\to Y$. By composition, this gives another plot $g\circ\alpha\in\diffeology_Z$. Now, since $f$ is a subduction, for every $t\in U_\alpha$ we can find a plot $\beta:V\to X$ defined on an open neighbourhood $t\in V\subseteq U_\alpha$ such that $g\circ\alpha|_V=f\circ \beta$. This gives a plot $(\beta,\alpha|_V):V\to X\times_Z Y$ that satisfies $\pr_2|_{X\times_ZY}\circ (\beta,\alpha|_V)=\alpha|_V$. The result follows by \cref{lemma:characterisation of subductions}.
		\end{proof}
	\end{lemma}
	
	The next result shows how two subductions interact with fibred products:
	\begin{lemma}
		\label{lemma:subduction and fibred product}
		Consider the following two commuting triangles of diffeological spaces and smooth maps:
		\begin{equation*}
			\begin{tikzcd}[column sep = small]
				{X_1} \arrow[rr, "f"] \arrow[rd, "r"'] &   & {Y_1} \arrow[ld, "R"] \\
				& A &                  
			\end{tikzcd}
			\qquad\text{and}\qquad
			\begin{tikzcd}[column sep = small]
				{X_2} \arrow[rr, "g"] \arrow[rd, "l"'] &   & {Y_2} \arrow[ld, "L"] \\
				& A, &                  
			\end{tikzcd}
		\end{equation*}
		where both $f$ and $g$ are subductions. Then the map
		\begin{equation*}
			(f\times g)|_{{X_1}\times_A{X_2}}:{X_1}\times_A^{r,l}{X_2}\longrightarrow {Y_1}\times_A^{R,L}{Y_2};\qquad (x_1,x_2)\longmapsto (f(x_1),g(x_2))
		\end{equation*}
		is also a subduction.
		
		\begin{proof}
			Clearly $f\times g$ is smooth, so we are left to show that the second condition in \cref{lemma:characterisation of subductions} is fulfilled. For that, take a plot $(\alpha_1,\alpha_2):U\to Y_1\times_A^{R,L}Y_2$, i.e., we have two plots $\alpha_1\in\diffeology_{Y_1}$ and $\alpha_2\in\diffeology_{Y_2}$ such that $R\circ \alpha_1=L\circ\alpha_2$. Now fix a point $t\in U$ in the domain. Then since both $f$ and $g$ are subductive, we can find two plots $\beta_1:U_1\to X_1$ and $\beta_2:U_2\to X_2$, defined on open neighbourhoods of $t\in U$, such that $\alpha_1|_{U_1}=f\circ \beta_1$ and $\alpha_2|_{U_2}=g\circ \beta_2$. Now the plot
			\begin{equation*}
				\left(
				\beta_1|_{U_1\cap U_2},\beta_2|_{U_1\cap U_2}
				\right)
				:
				U_1\cap U_2
				\longrightarrow 
				X_1\times X_2
			\end{equation*}
			takes values in the fibred product because
			\begin{equation*}
				r\circ \beta_1|_{U_2}
				=
				R\circ f\circ \beta_1|_{U_2}
				=
				R\circ \alpha_1|_{U_1\cap U_2}
				=
				L\circ \alpha_2|_{U_1\cap U_2}
				=
				l\circ \beta_2|_{U_1},
			\end{equation*}
			and we see that it lifts $(\alpha_1,\alpha_2)|_{U_1\cap U_2}$ along $f\times g$. 
		\end{proof}		
	\end{lemma}
	By setting $A=\lbrace\ast\rbrace$ to be the one-point space, this lemma gives in particular that the product $f\times g$ of two subductions is again a subduction.
	
	To end this section, we should also mention the existence of the notion of a \emph{local subduction} (also called \emph{strong subductions}):
	\begin{definition}
		A smooth surjection $f:X\to Y$ is called a \emph{local subduction} if for every \emph{pointed plot} of the form $\alpha:(U_\alpha,0)\to (Y,f(x))$ there exists a pointed plot $\beta:(V,0)\to (X,x)$, defined on an open neighbourhood $0\in V\subseteq U_\alpha$, such that $\alpha|_V=f\circ \beta$. 	
	\end{definition}
	Compare this to a definition of a subduction, where in general the plot $\beta$ does not have to hit the point $x$ in the domain of $f$. Note also that \emph{local} subduction does not mean \emph{locally a subduction everywhere}. 
	
	\begin{proposition}[{\cite[Article 2.16]{iglesias2013diffeology}}]
		\label{proposition:local subductions are surjective submersions}
		The local subductions between smooth manifolds are exactly the surjective submersions.
	\end{proposition}
	Due to the above proposition, the notion of a local subduction will be of interest when studying Lie groupoids in the framework of diffeological Morita equivalence we develop below. See \cref{section:diffeological bibundles between Lie groupoids}.
	
\section{Diffeological Groupoids}
	\label{section:diffeological groupoids}
	We assume that the reader is familiar with the definition of a (Lie) groupoid. A textbook reference for that theory is \cite{mackenzie2005general}. To fix our notation, we give here an informal description of a set-theoretic groupoid. A \emph{groupoid} consists of two sets: $G_0$ and $G$, together with five \emph{structure maps}. A groupoid will be denoted $\grpd{G}$, or just $G$. Here $G_0$ is the set of objects of the groupoid, and $G$ is the set of arrows. The five structure maps are
	\begin{enumerate}
		\item The \emph{source map} $\source:G\to G_0$, 
		\item The \emph{target map} $\target:G\to G_0$,
		\item The \emph{unit map} $\unit:G_0\to G$, mapping $x\mapsto \id_x$,
		\item The \emph{inversion map} $\inverse:G\to G$, mapping $g\mapsto g^{-1}$,
		\item And the \emph{composition}:
		\begin{equation*}
			\comp
			:
			G\times_{G_0}^{\source,\target}G
			\longrightarrow 
			G;
			\qquad
			(g,h)\mapsto g\circ h.
		\end{equation*}
	\end{enumerate}
	The composition is associative, and the identities and inverses behave as such. We say $\grpd{G}$ is a \emph{Lie groupoid} if both $G$ and $G_0$ are smooth manifolds such that the source and target maps are submersions, and each of the other structure maps are smooth. The definition of a diffeological groupoid is a straightforward generalisation of this:
	
	\begin{definition}
		\label{definition:diffeological groupoid}
		A \demph{diffeological groupoid} is a groupoid internal to the category of diffeological spaces. Concretely, this means that it is a groupoid $\grpd{G}$ such that the object space $G_0$ and arrow space $G$ are endowed with diffeologies that make all of the structure maps smooth. 
	\end{definition}
	
	As diffeology subsumes smooth manifolds, so do diffeological groupoids capture Lie groupoids. Note the main difference with the definition of a Lie groupoid is that we put no extra assumptions on the source and target maps. However:
	
	\begin{proposition}
		\label{proposition:source map is subduction}
		The source and target maps of a diffeological groupoid are subductions.
		
		\begin{proof}
			The smooth structure map $\unit:G_0\to G$, sending each object to its identity arrow, is a global smooth section of the source map, and hence by \cref{lemma:properties of subductions}\emph{(2)} the source map must be a subduction. Since the inversion map is a diffeomorphism, it follows that the target map is a subduction as well.
		\end{proof}
	\end{proposition}
	
	\begin{definition}
		\label{definition:isotropy groups}
		Let $\grpd{G}$ be a diffeological groupoid. The \demph{isotropy group} at $x\in G_0$ is the collection $G_x$ consisting of all arrows in $G$ from and to $x$:
		\begin{equation*}
			G_x:= \Hom_G(x,x)=\source^{-1}(\lbrace x\rbrace)\cap \target^{-1}(\lbrace x\rbrace).
		\end{equation*}
	\end{definition}
	
	\begin{definition}
		\label{definition:groupoid orbit space}
		Let $\grpd{G}$ be a diffeological groupoid. The \emph{orbit} of an object $x\in G_0$ is defined as
		\begin{equation*}
			\Orb_G(x)
			:=
			\lbrace y\in G_0 :\exists x\xrightarrow{~g~}y\rbrace
			=
			\target(\source^{-1}(\lbrace x\rbrace)).
		\end{equation*}
		The \demph{orbit space} of the groupoid is the space $G_0/G$ consisting of these orbits. We furnish the orbit space with the quotient diffeology from \cref{definition:quotient diffeology}, so that $\Orb_G:G_0\to G_0/G$ is a subduction.
	\end{definition}
	
	The orbit space of a Lie groupoid is not necessarily (canonically) a smooth manifold. The flexibility of diffeology allows us to study the smooth structure of orbit spaces of all diffeological groupoids. Below we give some examples of diffeological groupoids. 
	
	\begin{example}
		\label{example:relation groupoid}
		Let $X$ be a diffeological space, and let $R$ be an equivalence relation on $X$. We define the \demph{relation groupoid} $X\times_R X\rra X$ as follows. The space of arrows consists of exactly those pairs $(x,y)\in X\times X$ such that $xRy$. With the composition $(z,y)\circ(y,x):=(z,x)$, this becomes a diffeological groupoid. The orbit space $X/(X\times_R X)$ is just the quotient $X/R$. When $X$ is a smooth manifold, the relation groupoid becomes a Lie groupoid (even when the quotient is not a smooth manifold). 
	\end{example}
	
	\begin{example}
		\label{example:isotropy groupoid}
		Let $\grpd{G}$ be a diffeological groupoid. We can then consider the subgroupoid of $G$ that only consists of elements in isotropy groups:
		\begin{equation*}
			I_G
			:=
			\bigcup_{x\in G_0}G_x 
			\subseteq G.
		\end{equation*}
		This becomes a diffeological groupoid $I_G\rra G_0$ called the \demph{isotropy groupoid}. This has been studied in \cite[Example 2.1.9]{bos2007groupoids} in the context of Lie groupoids. Note that if $\grpd{G}$ is a Lie groupoid, then generally $I_G$ is not a submanifold of $G$, so the isotropy groupoid may no longer be a Lie groupoid. 
		
	\end{example}
	
	\begin{example}
		The \emph{thin fundamental groupoid} (or \emph{path groupoid}) $\Pi^\mathrm{thin}(M)$ of any smooth manifold $M$ is a diffeological groupoid \cite[Proposition A.25]{collier2016parallel}.
	\end{example}
	
	\begin{example}
		The \emph{groupoid of $\Sigma$-evolutions} of a Cauchy surface is a diffeological groupoid \cite[Section II.2.2]{glowacki2019}.
	\end{example}
	
	\begin{example}
		For any smooth surjection $\pi:X\to B$ between diffeological spaces, the fibres $X_b:=\pi^{-1}(\lbrace b\rbrace)$ get the subset diffeology from $X$. We then have a diffeological groupoid $\mathbf{G}(\pi)\rra B$ called the \emph{structure groupoid}, whose space of arrows is defined as
		\begin{equation*}
			\mathbf{G}(\pi):=\bigcup_{a,b\in B}\diff(X_a,X_b).
		\end{equation*}
		Structure groupoids play an important r\^ole in the theory of diffeological fibre bundles \cite[Chapter 8]{iglesias2013diffeology}. In general, they are too big to be Lie groupoids. They also generalise the notion of a \emph{frame groupoid} for a smooth vector bundle. Related to this, in \cite[Section 3.4]{schaaf2020diffeology-groupoids-and-ME} structure groupoids are used to define a notion of \emph{smooth linear representations} for diffeological groupoids. 
	\end{example}
	
	\begin{example}
		Given a diffeological space $X$, the \emph{germ groupoid} $\Germ(X)\rra X$ consists of all \emph{germs} of local diffeomorphisms on $X$. Even if $X$ itself is a smooth manifold, this is generally not a Lie groupoid. Germ groupoids are used in \cite{iglesias2018noncommutative,iglesias2020quasifolds}. A detailed construction of the diffeological structure of this groupoid appears in \cite[Section 6.1]{schaaf2020diffeology-groupoids-and-ME}.
	\end{example}
	
\section{Diffeological Groupoid Actions and -Bundles}
	\label{section:diffeological groupoid actions and bundles}
	In the following two sections we generalise the theory of Lie groupoid bibundles to the diffeological setting. The development we present here (as in \cite[Chapter IV]{schaaf2020diffeology-groupoids-and-ME}) is analogous to the development of the Lie version, save that we need to find a suitable replacement for the notion of a surjective submersion. Some of the proofs from the Lie theory can be performed almost \emph{verbatim} in our setting. These proofs already appear in the literature in various places: \cite{blohmann2008stacky,delHoyo2012Lie,landsman2001bicategories,moerdijk2005Poisson}, and also in the different setting of \cite{meyer2015groupoids}. We adopt many definitions and proofs from those sources, and point out how the diffeological theory subtly differs from the Lie theory. This difference mainly stems from the existence of quotients and fibred products of diffeological spaces, whereas in the Lie theory more care has to be taken. Ultimately, this extra care is what leads to a restricted Morita theorem for Lie groupoids, whereas the diffeological theorem is more general. In this section specifically we introduce \emph{diffeological groupoid actions} and \emph{-bundles}, two notions that form the ingredients for the main theory on bibundles.
	
	\subsection{Diffeological groupoid actions} 
	\label{section:diffeological groupoid actions}
	The most basic notion for the upcoming theory is that of a \emph{groupoid action}. For diffeological groupoids, the definition is the same as for Lie groupoids:
	
	\begin{definition}
		\label{definition:diffeological groupoid actions}
		Take a diffeological groupoid $\grpd{G}$, and a diffeological space $X$. A \emph{smooth left groupoid action} of $G$ on $X$ \emph{along} a smooth map $l_X:X\to G_0$ is a smooth function
		\begin{equation*}
			G\times_{G_0}^{\source,l_X}X
			\longrightarrow
			X;
			\qquad
			(g,x)
			\longmapsto
			g\cdot x,
		\end{equation*}
		satisfying the following three conditions:
		\begin{enumerate}
			\item For $g\in G$ and $x\in X$ such that $\source(g)=l_X(x)$ we have $l_X(g\cdot x)=\target(g)$.
			
			\item For every $x\in X$ we have $\id_{l_X(x)}\cdot x=x$.
			
			\item We have $h\cdot (g\cdot x)= (h\circ g)\cdot x$ whenever defined, i.e. when $\source(g)=l_X(x)$ and the arrows are composable.
		\end{enumerate}
		The smooth map $l_X:X\to G_0$ is called the \emph{left moment map}. In-line, we denote an action by $G\LAalong{l_X}X$. To save space, we may write $(g,x)\mapsto gx$ instead. 
		
		Right actions are defined similarly: a \emph{smooth right groupoid action} of $G$ on $X$ \emph{along} $r_X:X\to G_0$ is a smooth map
		\begin{equation*}
			X\times_{G_0}^{r_X,\target}G
			\longrightarrow
			X;
			\qquad
			(x,g)
			\longmapsto
			xg,
		\end{equation*}
		satisfying $r_X(xg)=\source(g)$, $x\cdot\id_{r_X(x)}=x$ and $(x\cdot g)\cdot h=x\cdot(g\circ h)$ whenever defined. Note how the r\^ole of the source and target maps are switched with respect to the definition of a left action. Right actions will be denoted by $X\RAalong{r_X}G$, and $r_X$ is called the \emph{right moment map}.
	\end{definition}
	
	\begin{example}
		\label{example:groupoid action on itself}
		Any diffeological groupoid $\grpd{G}$ acts on its own arrow space from the left and right by composition, which gives actions $G\LAalong{\target}G$ and $G\RAalong{\source}G$ that are both defined by $(g,h)\mapsto g\circ h$. 
	\end{example}
	\begin{definition}
		The \emph{orbit} of a point $x\in X$ in the space of an action $G\LAalong{l_X}X$ is defined as
		\begin{equation*}
			\Orb_G(x)
			:=
			\lbrace
			gx
			:
			g\in\source^{-1}(\lbrace l_X(x)\rbrace)
			\rbrace.
		\end{equation*}
		The \emph{quotient space} (or \emph{orbit space}) of the action is defined as the collection of all orbits, and denoted $X/G$. With the quotient diffeology, the \emph{orbit projection map} $\Orb_G:X\to X/G$ becomes a subduction.
	\end{definition}
	
	The following gives a notion of morphism between actions:
	\begin{definition}
		\label{definition:equivariant maps}
		Consider two smooth groupoid actions $G\LAalong{l_X}X$ and $G\LAalong{l_Y}Y$. A smooth map $\varphi:X\to Y$ is called \emph{$G$-equivariant} if $l_X=l_Y\circ\varphi$ and it commutes with the actions whenever defined: $\varphi(gx)=g\varphi(x)$. 
	\end{definition}
	
	\begin{definition}
		\label{definition:action category}
		The \emph{(smooth left) action category} $\Act(\grpd{G})$ of a diffeological groupoid $\grpd{G}$ is the category consisting of smooth left actions $G\LAalong{l_X}X$ as objects, and $G$-equivariant maps as morphisms. This forms the analogue of the category of (left) modules from ring theory. We show in \cref{section:invariance of representations} that the action category is in some sense a Morita invariant.
	\end{definition}
	
	\subsubsection{The balanced tensor product}
	We now give an important construction that will later allow us to define the \emph{composition} of bibundles.
	
	\begin{construction}
		\label{construction:balanced tensor product}
		Consider a diffeological groupoid $\grpd{H}$, with a smooth left action $H\LAalong{l_Y}Y$ and a smooth right action $X\RAalong{r_X}H$. On the fibred product $X\times_{H_0}^{r_X,l_Y}Y$ we define the following smooth right $H$-action. The moment map is $R:=r_X\circ \pr_1|_{X\times_{H_0}Y}=l_Y\circ\pr_2|_{X\times_{H_0}Y}$, and the action is given by:
		\begin{equation*}
			\left(
			X\times_{H_0}^{r_X,l_Y}Y
			\right)
			\times_{H_0}^{R,\target}H
			\longrightarrow
			X\times_{H_0}^{r_X,l_Y}Y;
			\qquad
			\left((x,y),h\right)
			\longmapsto
			(x\cdot h,h^{-1}\cdot y).
		\end{equation*}
		It is clear that this action is also smooth, and we call it the \emph{diagonal $H$-action}. The \emph{balanced tensor product} is the diffeological space defined as the orbit space of this smooth groupoid action:
		\begin{equation*}
			X\tensor_H Y 
			:= 
			\left(X\times_{H_0}^{r_X,l_Y}Y\right)/{H}.
		\end{equation*}
		The orbit of a pair of points $(x,y)$ in the balanced tensor product will be denoted $x\tensor y$. Whenever we encounter a term of the form $x\tensor y\in X\tensor_H Y$, we assume that it is well defined, i.e. $r_X(x)=l_Y(y)$. The terminology is explained by the following useful identity:
		\begin{equation*}
			xh\tensor y
			=
			x\tensor hy.
		\end{equation*}
		In the literature on Lie groupoids, this space is sometimes called the \emph{Hilsum-Skandalis tensor product}, named after a construction appearing in \cite{hilsum1987morphismes}.
	\end{construction}
	
	We note that this marks the first difference with the development of the Lie theory of bibundles and Morita equivalence. There, the balanced tensor product can only be defined when both $X\times_{H_0}^{r_X,l_Y}Y$ and the quotient by the diagonal $H$-action are smooth manifolds. This is usually only done after (bi)bundles are defined, and some principality conditions are presupposed. The principality then exactly ensures the existence of canonical differentiable structures on the fibred product and quotient. Here, the flexibility of diffeology allows us to define the balanced tensor product in an earlier stage of the development, and we do so to demonstrate this conceptual difference.

	\subsection{Diffeological groupoid bundles}
	\label{section:diffeological groupoid bundles}
	A groupoid bundle is a smooth map, whose domain carries a groupoid action, such that the fibres of the map are preserved by this action:
	
	\begin{definition}
		A \emph{smooth left diffeological groupoid bundle} is a smooth left groupoid action $G\LAalong{l_X}X$ together with a \emph{$G$-invariant} smooth map $\pi:X\to B$. We denote such bundles by $G\LAalong{l_X}X\xrightarrow{\pi}B$, and also call them \emph{(left) $G$-bundles}. \emph{Right} bundles are defined similarly, and denoted $B\xleftarrow{\pi}X\RAalong{r_X}G$. 
	\end{definition}
	
	The next definition gives a notion of morphism between bundles:
	
	\begin{definition}
		Consider two left $G$-bundles $G\LAalong{l_X}X\xrightarrow{\pi_X}B$ and $G\LAalong{l_Y}Y\xrightarrow{\pi_Y}B$ over the same base. A \emph{$G$-bundle morphism} is a $G$-equivariant smooth map $\varphi:X\to Y$ such that $\pi_X=\pi_Y\circ\varphi$. We make a similar definition for right bundles.
	\end{definition}
	
	In order to define Morita equivalence, we need to define a notion of when a bundle is \emph{principal}. For Lie groupoid bundles, these generalise the ordinary notion of smooth principal bundles of Lie groups and manifolds. That definition involves the notion of a surjective submersion. As we have mentioned, this notion needs to be generalised to diffeology. \cref{proposition:local subductions are surjective submersions} suggests that we could take \emph{local subductions}, since they directly generalise the surjective submersions. However, it turns out that \emph{subductions} behave sufficiently like submersions for the theory to work. The following definition then generalises the fact that the underlying bundle of a principal Lie groupoid bundle has to be a submersion:
	
	\begin{definition}
		A diffeological groupoid bundle $G\LAalong{l_X}X\xrightarrow{\pi}B$ is called \emph{subductive} if the underlying map $\pi:X\to B$ is a subduction.
	\end{definition}
	
	The following generalises the fact that the action of a principal Lie groupoid bundle has to be free and transitive on the fibres:
	
	\begin{definition}
		A diffeological groupoid bundle $G\LAalong{l_X}X\xrightarrow{\pi}B$ is called \emph{pre-principal} if the \emph{action map} $A_G:G\times_{G_0}^{\source,l_X}X\to X\times_B^{\pi,\pi}X$ mapping $(g,x)\mapsto (gx,x)$ is a diffeomorphism.
	\end{definition}
	
	Combining these two:
	\begin{definition}
		A diffeological groupoid bundle is called \emph{principal} if it is both subductive and pre-principal.
	\end{definition}
	
	This definition serves as our generalisation of principal Lie groupoid bundles, cf. \cite[Definition 2.10]{blohmann2008stacky} and \cite[Section 3.6]{delHoyo2012Lie}. Clearly any principal Lie groupoid bundle in the sense described in those references is also a principal diffeological groupoid bundle. Note that in the Lie theory, most constructions (such as the balanced tensor product) depend on the submersiveness of the underlying bundle map, so it makes little sense to define pre-principality for Lie groupoids. However, as we have already seen, in the diffeological case these constructions can be carried out more generally, and this will allow us to see what parts of the development of the theory depend on either the subductiveness or pre-principality of the bundles, rather than on full principality. In our development of the theory, some proofs can therefore be performed separately, whereas in the Lie theory they have to be performed at once. We hope this makes for clearer exposition.
	
	Note also that when a bundle $G\LAalong{l_X}X\xrightarrow{\pi}B$ is pre-principal, the action map induces a diffeomorphism $X/{\pi}\cong X/G$, and when the bundle is subductive, \cref{example:quotient by surjection or subduction} gives a diffeomorphism $B\cong X/{\pi}$. For a principal bundle we therefore have $B\cong X/G$. 
	
	\begin{example}
		The action of any diffeological groupoid $\grpd{G}$ on its own arrow space (\cref{example:groupoid action on itself}) forms a bundle $G\LAalong{\target}\xrightarrow{\source}G_0$. From \cref{proposition:source map is subduction} it follows that this bundle is principal.
	\end{example}
	
	\subsubsection{The division map of a pre-principal bundle}
	\label{section:division map}
	The material in this section is similar to \cite[Section 3.1]{blohmann2008stacky} for Lie groupoids. If a bundle $G\LAalong{l_X}X\xrightarrow{\pi}B$ is pre-principal, the fact that the action map is bijective gives that the action $G\LAalong{l_X}X$ has to be \emph{free}, and \emph{transitive} on the $\pi$-fibres. This means that for every two points $x,y\in X$ such that $\pi(x)=\pi(y)$, there exists a \emph{unique} arrow $g\in G$ such that $gy=x$. We denote this arrow by $\langle x,y\rangle_G$, and the map $\langle\cdot,\cdot\rangle_G$ is called the \emph{division map}\footnote{The notational resemblance to an inner-product is not accidental. The division map plays a very similar r\^ole to the inner product of a Hilbert \Cstar-module. For more on this analogy, see \cite[Section 3]{blohmann2008stacky}.}: 
	
	\begin{definition}
		Let $G\LAalong{l_X}X\xrightarrow{\pi}B$ be a pre-principal $G$-bundle, and let $A_G$ denote its action map. Then the \demph{division map} associated to this bundle is the smooth map
		\begin{equation*}
			\langle\cdot ,\cdot\rangle_G
			:
			X\times_B^{\pi,\pi}X
			\xrightarrow{\quad A_G^{-1}\quad } 
			G\times_{G_0}^{\source,l_X}X
			\xrightarrow{\quad \left.\pr_1\right|_{G\times_{G_0}X}\quad }
			G.
		\end{equation*}
	\end{definition}
	
	
	We summarise some algebraic properties of the division map that will be used in our proofs throughout later sections. The proofs are straightforward, and use the uniqueness property described above.
	
	\begin{proposition}
		\label{proposition:properties of division map}
		Let $G\LAalong{l_X}X\xrightarrow{\pi}B$ be a pre-principal $G$-bundle. Its division map $\langle\cdot,\cdot\rangle_G$ satisfies the following properties:
		\begin{enumerate}
			\item The source and targets are $\source(\langle x_1,x_2\rangle_G)=l_X(x_2)$ and $\target(\langle x_1,x_2\rangle_G)=l_X(x_1)$. 
			\item The inverses are given by $\langle x_1,x_2\rangle_G^{-1}=\langle x_2,x_1\rangle_G$.
			\item For every $x\in X$ we have $\langle x,x\rangle_G=\id_{l_X(x)}$.
			\item Whenever well-defined, we have $\langle gx_1,x_2\rangle_G=g\circ\langle x_1,x_2\rangle_G$.
		\end{enumerate}
	\end{proposition}
	
	\begin{proposition}
		Let $\varphi:X\to Y$ be a bundle morphism between two pre-principal $G$-bundles $G\LAalong{l_X}X\xrightarrow{\pi_X}B$ and $G\LAalong{l_Y}Y\xrightarrow{\pi_Y}B$. Denoting the division maps of these bundles respectively by $\langle\cdot,\cdot\rangle_G^X$ and $\langle\cdot,\cdot\rangle_G^Y$, we have for all $x_1,x_2\in X$ in the same $\pi_X$-fibre that:
		\begin{equation*}
			\langle x_1,x_2\rangle_G^X
			=
			\langle \varphi(x_1),\varphi(x_2)\rangle_G^Y.
		\end{equation*}
		
		\begin{proof}
			Note $\langle \varphi(x_1),\varphi(x_2)\rangle_G^Y$ is the unique arrow such that $\langle \varphi(x_1),\varphi(x_2)\rangle_G^Y\varphi(x_2)=\varphi(x_1)$. However, by $G$-equivariance we get $\varphi(x_1)=\varphi\left(\langle x_1,x_2\rangle_G^Xx_2\right)=\langle x_1,x_2\rangle_G^X\varphi(x_2)$, from which the claim immediately follows.
		\end{proof}
	\end{proposition}

	\subsubsection{Invertibility of $G$-bundle morphisms}
	\label{section:invertibility of bundle morphisms}
	We now prove a result that generalises the fact that morphisms between principal Lie group bundles are always diffeomorphisms. In our case we shall do the proof in two separate lemmas.
	
	\begin{lemma}
		\label{lemma:bundle morphism injective}
		Consider a $G$-bundle morphism $\varphi:X\to Y$ between a pre-principal bundle ${G\LAalong{l_X}X\xrightarrow{\pi_X}B}$ and a bundle $G\LAalong{l_Y}Y\xrightarrow{\pi_Y}B$ whose underlying action $G\LAalong{l_Y} Y$ is free. Then $\varphi$ is injective.
		
		\begin{proof}
			Since $G\LAalong{l_X} X\xrightarrow{\pi_X} B$ is pre-principal, we get a smooth division map $\langle \cdot,\cdot\rangle_G^X$. To start the proof, suppose that we have two points $x_1,x_2\in X$ satisfying $\varphi(x_1)=\varphi(x_2)$. Since $\varphi$ preserves the fibres, we get that
			\begin{equation*}
				\pi_X(x_1)=\pi_Y\circ\varphi(x_1)=\pi_Y\circ\varphi(x_2)=\pi_X(x_2).
			\end{equation*}
			Hence the pair $(x_1,x_2)$ defines an element in $X\times_B X$, so we get an arrow $\langle x_1,x_2\rangle_G^X\in G$, satisfying $\langle x_1,x_2\rangle_G^Xx_2=x_1$. If we apply $\varphi$ to this equation and use its $G$-equivariance, we get $\varphi(x_1)=\langle x_1,x_2\rangle_G^X\varphi(x_2)$. However, by assumption, $\varphi(x_1)=\varphi(x_2)$ and the action $G\LAalong{l_Y} Y$ is free, so we must have that $\langle x_1,x_2\rangle_G^X$ is the identity arrow at $l_Y\circ\varphi(x_2)=l_X(x_2)$. Hence we get the desired result:
			\begin{equation*}
				x_1=\langle x_1,x_2\rangle_G^X x_2=\id_{l_X(x_2)}x_2=x_2.
				\qedhere
			\end{equation*}
		\end{proof}
	\end{lemma}
	
	\begin{lemma}
		\label{lemma:bundle morphism subduction}
		Consider a $G$-bundle morphism $\varphi:X\to Y$ from a subductive bundle $G\LAalong{l_X}X\xrightarrow{\pi_X}B$ to a pre-principal bundle $G\LAalong{l_Y}Y\xrightarrow{\pi_Y}B$. Then $\varphi$ is a subduction.
		
		\begin{proof}
			Denote the smooth division map of $G\LAalong{l_Y} Y\xrightarrow{\pi_Y} B$ by $\langle \cdot,\cdot\rangle_G^Y$. Then $\varphi$ and $\langle \cdot,\cdot\rangle_G^Y$ combine into a smooth map
			\begin{equation*}
				\psi:
				X\times_B^{\pi_X,\pi_Y}Y\longrightarrow X;
				\qquad
				(x,y)\longmapsto \langle y,\varphi(x)\rangle_G^Yx.
			\end{equation*}
			Note that this is well-defined because if $\pi_X(x)=\pi_Y(y)$, then $\pi_Y\circ\varphi(x)=\pi_Y(y)$ as well, and moreover $l_Y\circ\varphi(x)=l_X(x)$, showing that the action on the right hand side is allowed. The $G$-equivariance of $\varphi$ then gives
			\begin{equation*}
				\varphi\circ \psi = \left.\pr_2\right|_{X\times_BY}.
			\end{equation*}
			Since $\pi_X$ is a subduction, so is $\pr_2|_{X\times_BY}$ by \cref{lemma:restriction of projection is subduction}, and by \cref{lemma:properties of subductions}\emph{(2)} it follows $\varphi$ is a subduction.
		\end{proof}
	\end{lemma}
	
	\begin{proposition}
		\label{proposition:bundle morphism on principal bundle is diffeomorphism}
		Any bundle morphism from a principal groupoid bundle to a pre-principal groupoid bundle is a diffeomorphism. In particular, both must then be principal.
		
		\begin{proof}
			By \cref{lemma:bundle morphism subduction} any such bundle morphism is a subduction, and since in particular the underlying action of a pre-principal bundle is free, it must also be injective by \cref{lemma:bundle morphism injective}. The result follows by \cref{proposition:injective subduction is diffeomorphism}. That the second bundle is principal too follows from the fact that a bundle map preserves the fibres, so the projection of the second bundle can be written as the composition of a diffeomorphism and a subduction.
		\end{proof}
	\end{proposition}

\section{Diffeological Bibundles and Morita Equivalence}
	\label{section:diffeological bibundles}
	This section contains the main definition of this paper: the notion of a \emph{biprincipal bibundle}, which immediately gives our definition of \emph{Morita equivalence}. The definition of groupoid bibundles for diffeology are a straightforward adaptation of the definition in the Lie case:
	
	\begin{definition}
		\label{definition:diffeological groupoid bibundles}
		Let $\grpd{G}$ and $\grpd{H}$ be two diffeological groupoids. A \emph{diffeological $(G,H)$-bibundle} consists of a smooth left action $G\LAalong{l_X}X$ and a smooth right action $X\RAalong{r_X}H$ such that the left moment map $l_X$ is $H$-invariant and the right moment map $r_X$ is $G$-invariant, and moreover such that the actions commute: $(g\cdot x)\cdot h=g\cdot (x\cdot h)$, whenever defined. We draw:
		\begin{equation*}
			\begin{tikzcd}[column sep = small]
				G\ar[r,phantom,"\laction"]\ar[d,shift left]\ar[d,shift right] & X\ar[dl,"l_X",pos=0.6]\ar[dr,"r_X", swap,pos=0.7] & \ar[l,phantom,"\raction"]H\ar[d,shift left]\ar[d,shift right]\\
				G_0 & & H_0,
			\end{tikzcd}
		\end{equation*}
		and denote them by $G\LAalong{l_X}X\RAalong{r_X}H$ in-line. Underlying each bibundle are two groupoid bundles: the \emph{left underlying $G$-bundle} $G\LAalong{l_X}X\xrightarrow{r_X}H_0$ and the \emph{right underlying $H$-bundle} $G_0\xleftarrow{l_X}X\RAalong{r_X}H$. It is the properties of these underlying bundles that will determine the behaviour of the bundle itself. 
	\end{definition}
	
	\begin{definition}
		\label{definition:left pre-principal}
		Consider a diffeological bibundle $G\LAalong{l_X}X\RAalong{r_X}H$. We say this bibundle is \demph{left pre-principal} if the left underlying bundle $G\LAalong{l_X}X\xrightarrow{r_X}H_0$ is pre-principal. We say it is \demph{right pre-principal} if the right underlying bundle $G_0\xleftarrow{l_X}X\RAalong{r_X}H$ is pre-principal. We make similar definitions for subductiveness and principality. Notice that, in this convention, if a bibundle $G\LAalong{l_X}X\RAalong{r_X}H$ is \emph{left} subductive, then its \emph{right} moment map $r_X$ is a subduction (and vice versa)\footnote{Note: \cite[Section 4.6]{delHoyo2012Lie} defines this differently, where ``[a] bundle is left (resp. right) principal if only the right (resp. left) underlying bundle is so.'' We suspect this may be a typo, since it apparently conflicts with their use of terminology in the proof of \cite[Theorem 4.6.3]{delHoyo2012Lie}. We stick to the terminology defined above, where \emph{left} principality pertains to the \emph{left} underlying bundle.}.
	\end{definition}
	
	We now have the main definition of this theory:
	\begin{definition}
		\label{definition:Morita equivalence and biprincipality}
		A diffeological bibundle is called: 
		\begin{enumerate}
			\item \demph{pre-biprincipal} if it is both left- and right pre-principal\footnote{The prefixes \emph{bi-} and \emph{pre-} commute: ``bi-(pre-principal)~=~pre-(biprincipal)''.};
			\item \demph{bisubductive} if it is both left- and right subductive;
			\item \demph{biprincipal} if it is both left- and right principal.
		\end{enumerate}
		Two diffeological groupoids $G$ and $H$ are called \demph{Morita equivalent} if there exists a biprincipal bibundle between them, and in that case we write $G\simeq_\mathrm{ME}H$.
	\end{definition}
	
	Compare this to the original definition \cite[Definition 2.1]{muhly1987equivalence} of equivalence for locally compact Hausdorff groupoids. We will prove in \cref{proposition:morita equivalence is equivalence relation} that Morita equivalence forms a genuine equivalence relation.
	\begin{example}
		\label{example:Lie ME is also diffeological ME}
		Since submersions between manifolds are subductions with respect to the manifold diffeologies, we see that if two \emph{Lie} groupoids $\grpd{G}$ and $\grpd{H}$ are Morita equivalent in the \emph{Lie} sense (e.g. \cite[Definition 2.15]{crainic2018orbispaces}), then they are Morita equivalent in the \emph{diffeological} sense. We remark on the converse question in \cref{section:diffeological bibundles between Lie groupoids}.
	\end{example}
	In fact, many elementary examples of Morita equivalences between Lie groupoids generalise straightforwardly to analogously defined diffeological groupoids. We refer to \cite[Section 4.3]{schaaf2020diffeology-groupoids-and-ME} for some of these examples. For us, the most important one is:
	
	\begin{example}
		\label{example:identity bibundle}
		Consider a diffeological groupoid $\grpd{G}$. There exists a canonical $(G,G)$-bibundle structure on the space of arrows $G$, which is called the \demph{identity bibundle}. The actions are just the composition in $G$ itself, as in \cref{example:groupoid action on itself}. Note that the identity bibundle is always biprincipal, because the action map has a smooth inverse $(g_1,g_2)\mapsto (g_1\circ g_2^{-1},g_2)$. This proves that any diffeological groupoid is Morita equivalent to itself, through the identity bibundle $G\LAalong{\target}G\RAalong{\source}G$. 
	\end{example}
	
	\begin{construction}
		\label{construction:opposite bibundle}
		Consider a diffeological bibundle $G\LAalong{l_X}X\RAalong{r_X}H$. The \emph{opposite bibundle} $H\LAalong{l_{\overline{X}}}\overline{X}\RAalong{r_{\overline{X}}}G$ is defined as follows. The underlying diffeological space does not change, $\overline{X}:=X$, but the moment maps switch, meaning that $l_{\overline{X}}:=r_X$ and $r_{\overline{X}}:=l_X$, and the actions are defined as follows:
		\begin{align*}
			H\LAalong{r_X}\overline{X};&\qquad h\cdot x:=xh^{-1},
			\\
			\overline{X}\RAalong{l_X}G;&\qquad x\cdot g:=g^{-1}x.
		\end{align*}
		Here the actions on the right-hand sides are the original actions of the bibundle. It is easy to see that performing this operation twice gives the original bibundle back. It is also important to note that for all properties defined in \cref{definition:left pre-principal}, taking the opposite merely switches the words `left' and `right'.
	\end{construction}
	The following extends \cref{proposition:properties of division map}\emph{(4)}:
	\begin{lemma}
		\label{lemma:opposite action division map}
		Consider a left pre-principal bibundle $G\LAalong{l_X}X\RAalong{r_X}H$, and also the opposite $G$-action $\overline{X}\RAalong{l_X}G$. Then, whenever defined, we have:
		\begin{equation*}
			\langle x_1,x_2 g\rangle_G = \langle x_1,x_2\rangle_G\circ g.
		\end{equation*}
		
		\begin{proof}
			This follows directly from \cref{proposition:properties of division map} and the definition of the opposite action:
			\begin{equation*}
				\langle x_1,x_2g\rangle_G
				=
				\langle x_1,g^{-1}x_2\rangle_G
				=
				\left(g^{-1}\circ\langle x_2,x_1\rangle_G\right)^{-1}
				=
				\langle x_1,x_2\rangle_G\circ g.\qedhere
			\end{equation*}
		\end{proof}
	\end{lemma}
	\subsection{Induced actions}
	\label{section:induced actions}
	A bibundle $G\laction X\raction H$ allows us to transfer a groupoid action $H\laction Y$ to a groupoid action $G\laction X\tensor_H Y$. This is called the \emph{induced action}, and, together with the balanced tensor product, will be crucial to define the composition of bibundles. The idea is that $G$ acts on the first component of $X\tensor_H Y$. 
	
	\begin{construction}
		\label{construction:induced action}
		Consider a diffeological bibundle $G\LAalong{l_X}X\RAalong{r_X}H$, and a smooth action $H\LAalong{l_Y}Y$. We construct a smooth left $G$-action on the balanced tensor product $X\tensor_H Y$. The left moment map is defined as
		\begin{equation*}
			L_X:X\tensor_H Y\longrightarrow G_0;
			\qquad 
			x\tensor y \longmapsto l_X(x).
		\end{equation*}
		This is well defined because $l_X$ is $H$-invariant, and smooth by \cref{lemma:properties of subductions}\emph{(3)}. For an arrow $g\in G$ with $\source(g)=L_X(x\tensor y)=l_X(x)$, define the action as:
		\begin{equation*}
			G\LAalong{L_X}X\tensor_H Y;\qquad g\cdot(x\tensor y):=(gx)\tensor y.
		\end{equation*}
		Note that the right hand side is well defined because $r_X$ is $G$-invariant, so $r_X(gx)=l_Y(y)$. Since there can be no confusion, we will drop all parentheses and write $gx\tensor y$ instead. That the action is smooth follows because $\left(g,(x,y)\right)\mapsto (gx,y)$ is smooth (on the appropriate domains) and by another application of \cref{lemma:properties of subductions}\emph{(3)}. Hence we obtain the \demph{induced action} $G\LAalong{L_X}X\tensor_H Y$.
		
		Now suppose that we are given a smooth $H$-equivariant map $\varphi:Y_1\to Y_2$ between two smooth actions $H\LAalong{l_1}Y_1$ and $H\LAalong{l_2}Y_2$. We define a map
		\begin{equation*}
			\id_X\tensor\varphi:X\tensor_HY_1\longrightarrow X\tensor_H Y_2;
			\qquad
			x\tensor y\longmapsto x\otimes\varphi(y).
		\end{equation*}
		The underlying map $X\times_{H_0}Y_1\to X\times_{H_0}Y_2:(x,y)\mapsto(x,\varphi(y))$ is clearly smooth. Then by composition of the projection onto $X\tensor_HY_2$ and \cref{lemma:properties of subductions}\emph{(3)}, we find $\id_X\tensor\varphi$ is smooth. Moreover, it is $G$-equivariant:
		\begin{equation*}
			\id_X\tensor\varphi(gx\tensor y)=gx\tensor\varphi(y)=g\left(\id_X\tensor\varphi(x\tensor y)\right).
		\end{equation*}
	\end{construction}
	
	\begin{definition}
		\label{definition:induced action functor}
		A diffeological bibundle $G\LAalong{l_X}X\RAalong{r_X}H$ defines an \demph{induced action functor}:
		\begin{align*}
			X\tensor_H-:\Act(\grpd{H})&\longrightarrow \Act(\grpd{G}),\\
			\left(H\LAalong{l_Y}Y\right)&\longmapsto \left(G\LAalong{L_X}X\tensor_H Y\right),\\
			\varphi&\longmapsto \id_X\tensor~\varphi.
		\end{align*}
		sending each smooth left $H$-action $\left(H\LAalong{l_Y}Y\right)\mapsto \left(G\LAalong{L_X}X\tensor_H Y\right)$ and each $H$-invariant map ${\varphi\mapsto \id_X\tensor\varphi}$. We will use this functor in \cref{section:invariance of representations}.
	\end{definition}

	\subsection{The bicategory of diffeological groupoids and -bibundles}
	Combining the balanced tensor product (\cref{construction:balanced tensor product}) and the induced action of a bibundle (\cref{construction:induced action}), we can define a notion of composition for diffeological bibundles, and thereby obtain a new sort of category of diffeological groupoids\footnote{The most straightforward way to obtain a (2-)category of diffeological groupoids is to consider the \emph{smooth functors} and \emph{smooth natural transformations}. We will not be studying this category in the current paper.}. Since performing multiple balanced tensor products is not strictly associative, we need to introduce a notion of comparison between diffeological bibundles.
	
	\begin{definition}
		\label{definition:biequivariant maps}
		Let $G\LAalong{l_X}X\RAalong{r_X}H$ and $G\LAalong{l_Y}Y\RAalong{r_Y}H$ be two bibundles between the same two diffeological groupoids. A smooth map $\varphi:X\to Y$ is called a \demph{bibundle morphism} if it is a bundle morphism between both underlying bundles. We also say that $\varphi$ is \demph{biequivariant}. Concretely, this means that the following diagram commutes:
		\begin{equation*}
			\begin{tikzcd}
				X \arrow[d, "l_X"'] \arrow[rd, "\varphi"] \arrow[r, "r_X"] & H_0                                  \\
				G_0                                                                    & Y, \arrow[l, "l_Y"] \arrow[u, "r_Y"']
			\end{tikzcd}
			\qquad
			\text{that is:}\qquad
			\begin{aligned}
				l_X&=l_Y\circ\varphi,\\
				r_X&=r_Y\circ\varphi,
			\end{aligned}
		\end{equation*}
		and that $\varphi$ is equivariant with respect to both actions. The isomorphisms of bibundles are exactly the diffeomorphic biequivariant maps. These are the 2-isomorphisms in $\DiffBiBund$. 
	\end{definition}
	
	The composition of bibundles is defined as follows:
	\begin{construction}
		\label{construction:bibundle composition}
		Consider two diffeological bibundles $G\LAalong{l_X}X\RAalong{r_X}H$ and $H\LAalong{l_Y}Y\RAalong{r_Y}K$. We shall define on $X\tensor_H Y$ a $(G,K)$-bibundle structure using the induced actions from \cref{construction:induced action}. 
		On the left we take the induced $G$-action along $L_X:X\tensor_H Y\to G_0$, which we recall maps ${x\tensor y \mapsto l_X(x)}$, defined by
		\begin{equation*}
			G\LAalong{L_X}X\tensor_H Y;\qquad g(x\tensor y):=(gx)\tensor y.
		\end{equation*}
		This action is well-defined because the $G$- and $H$-actions commute. 
		Similarly, we get an induced $K$-action on the right along $R_Y:X\tensor_H Y\to K_0$, which maps ${x\tensor y\mapsto r_Y(y)}$, given by
		\begin{equation*}
			X\tensor_H Y\RAalong{R_Y}K;\qquad (x\tensor y)k:= x\tensor(yk).
		\end{equation*}
		It is easy to see that these two actions form a bibundle $G\LAalong{L_X}X\tensor_HY\RAalong{R_Y}K$, which we also call the \demph{balanced tensor product}. Note that the moment maps are smooth by \cref{lemma:properties of subductions}\emph{(3)}.
	\end{construction}
	
	The following two propositions characterise the compositional structure of the balanced tensor product \emph{up to biequivariant diffeomorphism}. The first of these shows that the identity bibundle (\cref{example:identity bibundle}) is a \emph{weak identity:}
	\begin{proposition}
		\label{proposition:identity bibundle is weak identity}
		Let $G\LAalong{l_X}X\RAalong{r_X}H$ be a diffeological bibundle. Then there are biequivariant diffeomorphisms
		\begin{equation*}
			\begin{tikzcd}
				G\LAalong{L_G}G\tensor_G X\RAalong{R_X}H
				\arrow[d,Rightarrow,"\varphi"]
				\\
				G\LAalong{l_X}X\RAalong{r_X}H
			\end{tikzcd}
			\qquad\text{and}\qquad
			\begin{tikzcd}
				G\LAalong{L_X}X\tensor_H H\RAalong{R_H}H
				\arrow[d,Rightarrow,]
				\\
				G\LAalong{l_X}X\RAalong{r_X}H.
			\end{tikzcd}
		\end{equation*}
		
		\begin{proof}
			The idea of the proof is briefly sketched on \cite[p.8]{blohmann2008stacky}. The map $\varphi:G\tensor_G X\to X$ is defined by the action: $g\tensor x\mapsto gx$. This map is clearly well defined, and by an easy application of  \cref{lemma:properties of subductions}\emph{(3)} also smooth. Further note that $\varphi$ intertwines the left moment maps:
			\begin{equation*}
				l_X\circ\varphi(g\tensor x)=l_X(gx)=\target(g)=L_G(g\tensor x),
			\end{equation*}
			and similarly we find it intertwines the right moment maps. Associativity of the $G$-action and the fact that it commutes with the $H$-action directly ensure that $\varphi$ is biequivariant. Moreover, we claim that the smooth map $\psi:X\to G\tensor_G X$ defined by $x\mapsto \id_{l_X(x)}\tensor x$ is the inverse of $\varphi$. It follows easily that $\varphi\circ\psi=\id_X$, and the other side follows from the defining property of the balanced tensor product:
			\begin{equation*}
				\psi\circ\varphi(g\tensor x)=\psi(gx)=\id_{l_X(gx)}\tensor gx = (\id_{\target(g)}\circ g)\tensor x = g\tensor x.
			\end{equation*}
			It follows from an analogous argument that the identity bibundle of $H$ acts like a weak right inverse.
		\end{proof}
	\end{proposition}
	
	The second proposition shows that the balanced tensor is associative \emph{up to canonical biequivariant diffeomorphism:}
	\begin{proposition}
		\label{proposition:associativity of balanced tensor product}
		Let $G\LAalong{l_X}X\RAalong{r_X}H$, $H\LAalong{l_Y}Y\RAalong{r_Y}H'$, and $H'\LAalong{l_Z}Z\RAalong{r_Z}K$ be diffeological bibundles. Then there exists a biequivariant diffeomorphism
		\begin{equation*}
			\begin{tikzcd}
				G\LAalong{L_{X\tensor_HY}}\left(X\tensor_HY\right)\tensor_{H'}Z\RAalong{R_Z}K
				\arrow[d,Rightarrow,"A"]
				\\
				G\LAalong{L_X}X\tensor_H\left(Y\tensor_{H'}Z\right)\RAalong{R_{Y\tensor_{H'}Z}}K,
			\end{tikzcd}
			\quad
			A:(x\tensor y)\tensor z\longmapsto x\tensor(y\tensor z).
		\end{equation*}
		
		\begin{proof}
			That the map $A$ is smooth follows by \cref{lemma:properties of subductions}\emph{(3)}, because the corresponding underlying map ${\left((x,y),z\right)\mapsto \left(x,(y,z)\right)}$ is a diffeomorphism. The inverse of this diffeomorphism on the underlying fibred product induces exactly the smooth inverse of $A$, showing that $A$ is a diffeomorphism. Furthermore, it is easy to check that $A$ is biequivariant.
		\end{proof}
	\end{proposition}
	
	Combining \cref{proposition:identity bibundle is weak identity,proposition:associativity of balanced tensor product} gives that the balanced tensor product of bibundles behaves like the composition in a \emph{bicategory}. This is a category where the axioms of composition hold merely up to \emph{canonical 2-isomorphism}. For us, the 2-morphisms are the biequivariant smooth maps. For the precise definition of a bicategory we refer to e.g. \cite{macLane1998categories,lack2010companion}. The proof of the following is directly analogous to the one for the Lie theory, as explained throughout \cite{blohmann2008stacky}.
	
	\begin{theorem}
		\label{theorem:bicategory DiffBiBund}
		There is a bicategory $\DiffBiBund$ consisting of diffeological groupoids as objects, diffeological bibundles as morphisms with balanced tensor product as composition, and biequivariant smooth maps as 2-morphisms.
	\end{theorem}
	
	As we remarked in \cref{section:diffeological groupoid actions}, the balanced tensor product for Lie groupoids can only be constructed for \emph{left} (or \emph{right}) \emph{principal} bibundles. This means that in the Lie theory, the category of bibundles only consists of the left (or right) principal bibundles, since otherwise the composition cannot be defined. For diffeology we obtain a bicategory of \emph{all} bibundles.

	\subsection{Properties of bibundles under composition and isomorphism}
	\label{section:properties of bibundles under composition and isomorphism}
	We study how the properties of diffeological bibundles defined in \cref{definition:left pre-principal} are preserved under the balanced tensor product and biequivariant diffeomorphism. These results will be crucial in characterising the weakly invertible bibundles. First we show that left subductive and left pre-principal bibundles are closed under composition.
	
	\begin{proposition}
		\label{proposition:subductive balanced tensor product}
		The balanced tensor product preserves left subductiveness.
		
		\begin{proof}
			Consider the balanced tensor product $G\LAalong{L_X}X\tensor_H Y\RAalong{R_Y}K$ of two left subductive bibundles $G\LAalong{l_X}X\RAalong{r_X}H$ and $H\LAalong{l_Y}Y\RAalong{r_Y}K$. We need to show that the right moment map $R_Y:X\tensor_HY\to K_0$ is a subduction. But, note that it fits into the following commutative diagram:
			\begin{equation*}
				\begin{tikzcd}
					{X\times_{H_0}^{r_X,l_Y}Y} \arrow[r, "\pi"] \arrow[d, "\pr_2|_{X\times_{H_0}Y}"'] & X\tensor_HY \arrow[d, "R_Y"] \\
					Y \arrow[r, "r_Y"']                                             & K_0 .                        
				\end{tikzcd}
			\end{equation*}
			Here $\pi$ is the canonical quotient projection. The restricted projection $\pr_2|_{X\times_{H_0}Y}$ is a subduction by \cref{lemma:restriction of projection is subduction}, since $r_X$ is a subduction. Moreover, $r_Y$ is a subduction, so the bottom part of the diagram is a subduction. It follows by \cref{lemma:properties of subductions}\emph{(3)} that $R_Y$ is a subduction. 
		\end{proof}
	\end{proposition}
	
	Note that, even though $R_Y$ only explicitly depends on the moment map $r_Y$, the proof still depends on the subductiveness of $r_X$ as well.
	
	To prove that the balanced tensor product of two left pre-principal bibundles is again left pre-principal, we need the following lemma, describing how the division map interacts with the bibundle structure, extending the list in \cref{proposition:properties of division map} on the algebraic properties of the division map.
	
	\begin{lemma}
		\label{lemma:division map on bibundle}
		Let $G\LAalong{l_X}X\RAalong{r_X}H$ be a left pre-principal bibundle, and denote its division map by $\langle\cdot,\cdot\rangle_G$. Then, whenever defined:
		\begin{equation*}
			\langle x_1,x_2h\rangle_G=\langle x_1h^{-1},x_2\rangle_G,
			\qquad
			\text{or equivalently:}
			\qquad
			\langle x_1h,x_2h\rangle_G=\langle x_1,x_2\rangle_G.
		\end{equation*}
		
		\begin{proof}
			The arrow $\langle x_1h,x_2h\rangle_G\in G$ is the unique one so that $\langle x_1h,x_2h\rangle_G(x_2h)=x_1h$. Now, since the actions commute, we can multiply both sides of this equation from the right by $h^{-1}$, which gives $\langle x_1h,x_2h\rangle_X x_2=x_1$, and this immediately gives our result.
		\end{proof}
	\end{lemma}
	
	\begin{proposition}
		\label{proposition:pre-principal balanced tensor product}
		The balanced tensor product preserves left pre-principality.
		
		\begin{proof}
			To start the proof, take two left pre-principal bibundles, with our usual notation: $G\LAalong{l_X}X\RAalong{r_X}H$ and $H\LAalong{l_Y}Y\RAalong{r_Y}K$. Denote their division maps by $\langle\cdot,\cdot\rangle_G^X$ and $\langle \cdot,\cdot\rangle_H^Y$, respectively. Using these, we will construct a smooth inverse of the action map of the balanced tensor product. Let us denote the action map of the balanced tensor product by
			\begin{equation*}
				\Phi:G\times_{G_0}^{\source,L_X}\left(X\tensor_HY\right)
				\longrightarrow
				\left(X\tensor_HY\right)\times_{K_0}^{R_Y,R_Y}\left(X\tensor_HY\right),
			\end{equation*}
			mapping $(g,x\tensor y)\mapsto (gx\tensor y,x\tensor y)$. After some calculations (which we describe below), we propose the following map as an inverse for $\Phi$:
			\begin{align*}
				\Psi:
				\left(X\tensor_HY\right)\times_{K_0}^{R_Y,R_Y}\left(X\tensor_HY\right)
				&\longrightarrow
				G\times_{G_0}^{\source,L_X}\left(X\tensor_HY\right);
				\\
				\left(x_1\tensor y_1,x_2\tensor y_2\right)
				&\longmapsto
				\left(\left\langle x_1\langle y_1,y_2\rangle_H^Y,x_2 \right\rangle_G^X,x_2\tensor y_2
				\right).
			\end{align*}
			It is straightforward to check that every action and division occurring in this expression is well defined. We need to check that $\Psi$ is independent on the representations of $x_1\tensor y_1$ and $x_2\tensor y_2$. Only the first component $\Psi_1$ of $\Psi$ could be dependent on the representations, so we focus there. Suppose we have two arrows $h_1,h_2\in H$ satisfying $\target(h_i)=r_X(x_i)=l_Y(y_i)$, so that $x_ih_i\tensor h_i^{-1}y_i=x_i\tensor y_i$. For the division of $y_2$ and $y_1$ we then use \cref{proposition:properties of division map} to get:
			\begin{equation*}
				\langle h_1^{-1}y_1,h_2^{-1}y_2\rangle_H^Y
				=
				h_1^{-1}\circ\langle y_1,h_2^{-1}y_2\rangle_H^Y
				=
				h_1^{-1}\circ \left( h_2^{-1}\circ\langle y_2,y_1\rangle_H^Y\right)^{-1}
				=
				h_1^{-1}\circ \langle y_1,y_2\rangle_H^Y\circ h_2.
			\end{equation*}
			Then, using this and \cref{lemma:division map on bibundle}, we get:
			\begin{align*}
				\Psi_1(x_1h_1\tensor h_1^{-1}y_1,x_2h_2\tensor h_2^{-1}y_2)
				&=
				\left\langle x_1h_1\langle h_1^{-1}y_1,h_2^{-1}y_2\rangle_H^Y,x_2h_2\right\rangle_G^X
				\\&=
				\left\langle (x_1h_1)\left(h_1^{-1}\circ\langle y_1,y_2\rangle_H^Y\circ h_2\right),x_2h_2\right\rangle_G^X
				\\&=
				\left\langle (x_1\langle y_1,y_2\rangle)h_2,x_2h_2\right\rangle_G^X
				\\&=
				\left\langle x_1\langle y_1,y_2\rangle_H^Y,x_2\right\rangle_G^X.
			\end{align*}
			Since the second component of $\Psi$ is by construction independent on the representation, it follows that $\Psi$ is a well-defined function. We now need to show that $\Psi$ is smooth. The second component is clearly smooth, because it is just the projection onto the second component of the fibred product. That the other component is smooth follows from \cref{lemma:properties of subductions,lemma:subduction and fibred product}.
			Writing
			\begin{equation*}
				\psi:\left((x_1,y_1),(x_2,y_2)\right)\longmapsto \langle x_1\langle y_1,y_2\rangle_H^Y,x_2\rangle_G^X
			\end{equation*}			
			and $\pi:X\times_{H_0}^{r_X,l_Y}Y\to X\tensor_H Y$ for the canonical projection, we get a commutative diagram
			\begin{equation*}
				\begin{tikzcd}
					\left(X\times_{H_0}^{r_X,l_Y}Y\right)\times_{K_0}^{\overline{r_Y},\overline{r_Y}}\left(X\times_{H_0}^{r_X,l_Y}Y\right) \arrow[rr, "(\pi\times\pi)|_{\mathrm{\dom(\psi)}}"] \arrow[rd, "\psi"', bend right=10] &   & \left(X\tensor_HY\right)\times_{K_0}^{R_Y,R_Y}\left(X\tensor_HY\right) \arrow[ld, "\Psi_1",bend left=10] \\
					& G. &                          
				\end{tikzcd}
			\end{equation*}
			Here we temporarily use the notation $\overline{r_Y}:=r_Y\circ\pr_2|_{X\times_{H_0}Y}$, which satisfies $R_Y\circ\pi = \overline{r_Y}$. Therefore by \cref{lemma:subduction and fibred product} the top arrow in this diagram is a subduction. Since the map $\psi$ is evidently smooth, it follows by \cref{lemma:properties of subductions}\emph{(3)} that the first component $\Psi_1$, and hence $\Psi$ itself, must be smooth.
			
			Thus, we are left to show that $\Psi$ is an inverse for $\Phi$. That $\Psi$ is a right inverse for $\Phi$ now follows by simple calculation using \cref{proposition:properties of division map,lemma:division map on bibundle}:
			
			\begin{equation*}
				\Psi\circ\Phi(g,x\tensor y)
				=
				\Psi(gx\tensor y,x\tensor y)
				=
				\left(\langle gx\langle y,y\rangle_H^Y,x\rangle_G^X,x\tensor y\right)
				= 
				\left(g\circ \langle x,x\rangle_G^X,x\tensor y\right)
				=
				(g,x\tensor y).
			\end{equation*}
			For the other direction, we calculate:
			\begin{align*}
				\Phi\circ\Psi(x_1\tensor y_1,x_2\tensor y_2)
				&=
				\Phi\left(\left\langle x_1\langle y_1,y_2\rangle_H^Y,x_2 \right\rangle_G^X,x_2\tensor y_2
				\right)
				\\&=
				\left(\left\langle x_1\langle y_1,y_2\rangle_H^Y,x_2 \right\rangle_G^Xx_2\tensor y_2,x_2\tensor y_2\right)
				\\&=
				\left(x_1\langle y_1,y_2\rangle_H^Y\tensor y_2,x_2\tensor y_2\right)
				\\&=
				\left(x_1\tensor \langle y_1,y_2\rangle_H^Yy_2,x_2\tensor y_2\right)
				\\&=
				\left(x_1\tensor y_1,x_2\tensor y_2\right).
			\end{align*}
			Here in the second to last step we use the properties of the balanced tensor product to move the arrow $\langle y_1,y_2\rangle_H^Y$ over the tensor symbol. Hence we conclude that $\Phi$ is a diffeomorphism, which proves that $G\LAalong{L_X}X\tensor_H Y\RAalong{R_Y}K$ is a left pre-principal bibundle.
		\end{proof}
	\end{proposition}
	
	Next we show that left subductiveness and left pre-principality are also preserved under biequivariant diffeomorphism.
	
	\begin{proposition}
		\label{proposition:pre-principality and isomorphism}
		Left pre-principality is preserved by biequivariant diffeomorphism.
		
		\begin{proof}
			Suppose that $\varphi:X\to Y$ is a biequivariant diffeomorphism from a left pre-principal bibundle $G\LAalong{l_X}X\RAalong{r_X}H$ to another diffeological bibundle $G\LAalong{l_Y}Y\RAalong{r_Y}H$. Denote their left action maps by $A_X$ and $A_Y$, respectively. The following square commutes because of biequivariance:
			\begin{equation*}
				\begin{tikzcd}
					{G\times_{G_0}^{\source,l_X}X} \arrow[d, "(\id_G\times\varphi)|_{G\times_{G_0}X}"'] \arrow[r, "A_X"] & {X\times_{H_0}^{r_X,r_X}X} \arrow[d, "(\varphi\times\varphi)|_{X\times_{H_0}X}"] \\
					{G\times_{G_0}^{\source,l_Y}Y} \arrow[r, "A_Y"']                                                     & {Y\times_{H_0}^{r_Y,r_Y}Y.}                                                      
				\end{tikzcd}
			\end{equation*}
			It is easy to see that both vertical maps are diffeomorphisms. Hence it follows $A_Y$ must be a diffeomorphism as well.
		\end{proof}
	\end{proposition}
	
	\begin{proposition}
		\label{proposition:subductiveness and isomorphism}
		Left subductiveness is preserved by biequivariant diffeomorphism.
		
		\begin{proof}
			Suppose that $\varphi:X\to Y$ is a biequivariant diffeomorphism from a left subductive bibundle $G\LAalong{l_X}X\RAalong{r_X}H$ to $G\LAalong{l_Y}Y\RAalong{r_Y}H$. That the first bundle is left subductive means that $r_X$ is a subduction, but since $\varphi$ intertwines the moment maps, it follows immediately that $r_Y=r_X\circ\varphi^{-1}$ is a subduction as well.
		\end{proof}
	\end{proposition}
	
	Of course, these four propositions all hold for their respective `right' versions as well. This can be proved formally, without repeating the work, by using opposite bibundles.
	
	\begin{corollary}
		\label{proposition:morita equivalence is equivalence relation}
		Morita equivalence defines an equivalence relation between diffeological groupoids.
		
		\begin{proof}
			Morita equivalence is reflexive by the existence of identity bibundles, which are always biprincipal (\cref{example:identity bibundle}). It is also easy to check that the opposite bibundle (\cref{construction:opposite bibundle}) of a biprincipal bibundle is again biprincipal, showing that Morita equivalence is symmetric. Transitivity follows directly from \cref{proposition:pre-principal balanced tensor product,proposition:subductive balanced tensor product} and their opposite versions.
		\end{proof}
	\end{corollary}

	\subsection{Weak invertibility of diffeological bibundles}
	\label{section:weak invertibility of diffeological bibundles}
	In this section we prove the main Morita \cref{theorem:weakly invertible bibundles are the biprincipal ones}. As we explained in the \nameref{section:introduction}, in the bicategory of diffeological groupoids we get a notion of \emph{weak isomorphism}. Let us describe these explicitly: a bibundle $G\laction X\raction H$ is weakly invertible if and only if there exists a second bibundle $H\laction Y\raction G$, such that $X\tensor_H Y$ is biequivariantly diffeomorphic to $G$ and $Y\tensor_G X$ is biequivariantly diffeomorphic to $H$. The Morita theorem says that such a weak inverse exists if and only if the bibundle is biprincipal. Let us recall the corresponding statement in the Lie theory: a (say) left principal bibundle has a left principal weak inverse if and only if it is biprincipal \cite[Proposition 4.21]{landsman2001quantized}. Here both the original bibundle and its weak inverse have to be left principal, since everything takes place in a bicategory of Lie groupoids and left principal bibundles. According to \cref{theorem:bicategory DiffBiBund} we get a bicategory of arbitrary bibundles, and the question of weak invertibility becomes a slightly more general one, since we do not start out with a bibundle that is already left principal. Instead we have to infer left principality from bare weak invertibility, where neither the weak inverse may be assumed to be left principal.  
	
	One direction of the claim in the main theorem is relatively straightforward, and is the same as for Lie groupoids:
	
	\begin{proposition}
		\label{proposition:biprincipal bibundle is weakly invertible}
		Let $G\LAalong{l_X}X\RAalong{r_X}H$ be a biprincipal bibundle. Then its opposite bundle $H\LAalong{r_X}\overline{X}\RAalong{l_X}G$ is a weak inverse.
		
		\begin{proof}
			We construct biequivariant diffeomorphisms
			\begin{equation*}
				\begin{tikzcd}
					G\LAalong{L_X}X\tensor_H \overline{X}\RAalong{R_{\overline{X}}}G
					\arrow[d,Rightarrow,"\varphi_G"']
					\\
					G\LAalong{\target}G\RAalong{\source}G,
				\end{tikzcd}
				\quad\text{and}\quad
				\begin{tikzcd}
					H\LAalong{L_{\overline{X}}}\overline{X}\tensor_G X\RAalong{R_X}H
					\arrow[d,Rightarrow,"\varphi_H"']
					\\
					H\LAalong{\target}H\RAalong{\source}H.
				\end{tikzcd}
			\end{equation*}
			Since the original bundle is pre-biprincipal, we have a division map $\langle\cdot,\cdot\rangle_G:X\times_{H_0}^{r_X,r_X}\overline{X}\to G$. We define a new function
			\begin{equation*}
				\varphi_G
				:
				X\tensor_H\overline{X}
				\longrightarrow 
				G;
				\qquad
				x_1\tensor x_2
				\longmapsto 
				\langle x_1,x_2\rangle_G.
			\end{equation*}
			This is independent on the representation of the tensor product by \cref{lemma:division map on bibundle}, and smooth by \cref{lemma:properties of subductions}\emph{(3)} since $\varphi_G\circ\pi=\langle\cdot,\cdot\rangle_G$, where $\pi$ is the canonical projection onto the balanced tensor product. We check that $\varphi_G$ is biequivariant. It is easy to check that $\varphi_G$ intertwines the moment maps, for example:
			\begin{equation*}
				\source\circ\varphi_G(x_1\tensor x_2)
				=
				\source\left(\langle x_1,x_2\rangle_G\right)
				=
				l_X(x_2)
				=
				R_{\overline{X}}(x_1\tensor x_2).
			\end{equation*}
			The left $G$-equivariance of $\varphi_G$ follows directly out of \cref{proposition:properties of division map}, and the right $G$-equivariance follows from \cref{lemma:opposite action division map}. Hence $\varphi_G$ is a genuine bibundle morphism.
			
			Since the original bundle is biprincipal, so is its opposite, and therefore by \cref{proposition:pre-principal balanced tensor product,proposition:subductive balanced tensor product} it follows that both balanced tensor products are also biprincipal. Therefore $\varphi_G$ is in particular a left $G$-equivariant bundle morphism from a principal bundle $G\LAalong{L_X}X\tensor_H\overline{X}\xrightarrow{R_{\overline{X}}}G_0$ to a pre-principal bundle $G\LAalong{\target}G\xrightarrow{\source}G_0$, and hence a diffeomorphism by \cref{proposition:bundle morphism on principal bundle is diffeomorphism}. This proves that the opposite bibundle is a weak right inverse. Note that we already need full biprincipality of the original bibundle for this. To prove that it is also a weak left inverse we make an analogous construction for $\varphi_H$, which we leave to the reader.
		\end{proof}
	\end{proposition}
	
	The rest of this section will be dedicated to proving the converse of this claim, i.e., that a weakly invertible bibundle is biprincipal. First let us remark that by imitating a result from the Lie theory, we can obtain a partial result in this direction. Let us denote by $\DiffBiBund_\mathrm{LP}$ the bicategory of diffeological groupoids and left principal bibundles. Note that by \cref{section:properties of bibundles under composition and isomorphism} left principality is preserved by the balanced tensor product, so this indeed forms a subcategory. 
	\begin{theorem}
		\label{theorem:analogue}
		A left principal diffeological bibundle has a left principal weak inverse if and only if it is biprincipal. That is, the weakly invertible bibundles in $\DiffBiBund_\mathrm{LP}$ are exactly the biprincipal ones.
		
		\begin{proof}
			This follows by combining \cref{proposition:biprincipal bibundle is weakly invertible} with an adaptation of an argument from the Lie groupoid theory as in \cite[Proposition 2.9]{moerdijk2005Poisson}. A more detailed proof (for diffeological groupoids) is in \cite[Proposition 4.61]{schaaf2020diffeology-groupoids-and-ME}.
		\end{proof}
	\end{theorem}
	
	This theorem is the most direct analogue of \cite[Proposition 4.21]{landsman2001quantized} in the setting of diffeology. Our main theorem will be a further generalisation of this, which says that the same claim holds in the larger bicategory $\DiffBiBund$ of \emph{all} bibundles. We break the proof down in several steps, starting with the implication of bisubductiveness:
	
	\begin{proposition}
		\label{proposition:weakly invertible implies bisubductive}
		A weakly invertible diffeological bibundle is bisubductive.
		
		\begin{proof}
			Suppose we have a bibundle $G\LAalong{l_X}X\RAalong{r_X}H$ that admits a weak inverse $H\LAalong{l_Y}Y\RAalong{r_Y}G$. Let us denote the included biequivariant diffeomorphisms by $\varphi_G:X\tensor_H Y\to G$ and $\varphi_H:Y\tensor_G X\to H$, as usual. Since the identity bibundles of $G$ and $H$ are both biprincipal, it follows by \cref{proposition:subductiveness and isomorphism} that the moment maps $L_X$, $R_X$, $L_Y$ and $R_Y$ are all subductions. Together with the original moment maps, we get four commutative squares, each of the form:
			\begin{equation*}
				\begin{tikzcd}
					{X\times_{H_0}^{r_X,l_Y}Y} \arrow[d, "\pr_1|_{X\times_{H_0}Y}"'] \arrow[r, "\pi"] & X\tensor_H Y \arrow[d, "L_X"] \\
					X \arrow[r, "l_X"']                                                               & G_0.                          
				\end{tikzcd}
			\end{equation*}
			Here $\pi:X\times_{H_0}^{r_X,l_Y}Y\to X\tensor_HY$ is the quotient map of the diagonal $H$-action. By \cref{lemma:properties of subductions}\emph{(3)} it follows that, since $L_X$ is a subduction, so is $l_X\circ\pr_1|_{X\times_{H_0}Y}$, and in turn by \cref{lemma:properties of subductions}\emph{(2)} it follows $l_X$ is a subduction. In a similar fashion we find that $r_X$, $l_Y$ and $r_Y$ are all subductions as well.
		\end{proof}
	\end{proposition}

	This proposition gets us halfway to proving that weakly invertible bibundles are biprincipal. To prove that they are pre-biprincipal, it is enough to construct smooth division maps. We will give this construction below (\cref{construction:local division map}), which follows from a careful reverse engineering of the division map of a pre-principal bundle. Recall from \cref{proposition:pre-principal balanced tensor product} that the smooth inverse of the action map contains the information of both the $G$-division map and the $H$-division map. Specifically, the first component of the inverse is of the form $\langle x_1\langle y_1,y_2\rangle_H^Y,x_2\rangle_G^X$, in which if we set $y_1=y_2$, we simply reobtain the $G$-division map $\langle x_1,x_2\rangle_G^X$. The question is if this ``reobtaining'' can be done in a smooth way. This is not so obvious at first. Namely, if we vary $(x_1,x_2)$ smoothly within $X\times_{H_0}^{r_X,r_X}X$, can we guarantee that $y_1$ and $y_2$ vary smoothly with it, while still retaining the equalities $r_X(x_i)=l_Y(y_i)$ and $y_1=y_2$? The elaborate \cref{construction:local division map} proves that this can indeed be done. An essential part of our argument will be supplied by the following two lemmas.
	
	\begin{lemma}
		\label{lemma:actions of weakly invertible bibundle are free}
		When $G\LAalong{l_X}X\RAalong{r_X}H$ is a weakly invertible bibundle, admitting a weak inverse $H\LAalong{l_Y}Y\RAalong{r_Y}G$, then all four actions are free.
		
		\begin{proof}
			This follows from an argument that is used in the proof of \cite[Proposition 3.23]{blohmann2008stacky}. Suppose we have an arrow $h\in H$ and a point $y\in Y$ such that $hy=y$. By \cref{proposition:weakly invertible implies bisubductive} it follows that in particular $l_X$ is surjective, so we can find $x\in X$ such that $y\tensor x\in Y\tensor_G X$. Then 
			\begin{equation*}
				h(y\tensor x)=(hy)\tensor x=y\tensor x.
			\end{equation*}
			But by \cref{proposition:pre-principality and isomorphism} the bundle $H\LAalong{L_Y}Y\tensor_G X\xrightarrow{R_X}G_0$, which is equivariantly diffeomorphic to the identity bundle on $H$, is pre-principal. So, the left action $H\laction Y\tensor_G X$ is free, and hence $h=\id_{L_Y(y\tensor x)}=\id_{l_Y(y)}$, proving that $H\laction Y$ is also free. That the three other actions are free follows analogously.
		\end{proof}
	\end{lemma}
	
	\begin{lemma}
		\label{lemma:free action and balanced tensor product}
		Let $X\RAalong{r_X}H$ and $H\LAalong{l_Y}Y$ be smooth actions, so that we can form the balanced tensor product $X\tensor_H Y$. Suppose that $H\laction Y$ is free. Then $x_1\tensor y= x_2\tensor y$ if and only if $x_1=x_2$. Similarly, if $X\raction H$ is free, then $x\tensor y_1=x\tensor y_2$ if and only if $y_1=y_2$.
		
		\begin{proof}
			If $x_1=x_2$ to begin with, the implication is trivial. Suppose therefore that $x_1\tensor y=x_2\tensor y$, which means that there exists an arrow $h\in H$ such that $(x_1h^{-1},hy)=(x_2,y)$. In particular $hy=y$, which, because the action on $Y$ is free, implies $h=\id_{l_Y(y)}$, and it follows that $x_1=x_1 \id_{l_Y(y)}^{-1}=x_2$.
		\end{proof}
	\end{lemma}
	
	We shall now describe how the division map arises from local data:

	\begin{construction}
		\label{construction:local division map}
		For this construction to work, we start with a diffeological bibundle $G\LAalong{l_X}X\RAalong{r_X}H$, admitting a weak inverse $H\LAalong{l_Y}Y\RAalong{r_Y}G$. Then consider a pointed plot $\alpha:(U_\alpha,0)\to (X\times_{H_0}^{r_X,r_X}X,(x_1,x_2))$. We split $\alpha$ into the components $(\alpha_1,\alpha_2)$, which in turn are pointed plots $\alpha_i:(U_\alpha,0)\to (X,x_i)$ satisfying $r_X\circ \alpha_1=r_X\circ\alpha_2:U_\alpha\to H_0$. This equation gives a plot of $H_0$, and since by \cref{proposition:weakly invertible implies bisubductive} the moment map $l_Y:Y\to H_0$ is a subduction, for every $t\in U_\alpha$ we can find a plot $\beta:V\to Y$, defined on an open neighbourhood $t\in V\subseteq U_\alpha$, such that $r_X\circ\alpha_i|_V=l_Y\circ\beta$. From this equation it follows that the smooth maps $(\alpha_i|_V,\beta):V\to X\times_{H_0}^{r_X,l_Y}Y$ define two plots of the underlying space of the balanced tensor product. Applying the quotient map $\pi:X\times_{H_0}^{r_X,l_Y}Y\to X\tensor_H Y$, we thus get two full-fledged plots $s\mapsto \alpha_i|_V(s)\tensor\beta(s)$ of the balanced tensor product. We combine these two plots to define yet another smooth map:
		\begin{equation*}
			\Omega^\alpha|_V:=\left(\pi\circ\left(\alpha_1|_V,\beta\right),\pi\circ\left(\alpha_2|_V,\beta\right)\right)
			:
			V
			\longrightarrow
			\left(X\tensor_H Y\right)\times_{G_0}^{R_Y,R_Y}\left(X\tensor_H Y\right).
		\end{equation*}
		Note that $\Omega^\alpha|_V$ lands in the right codomain because $R_Y\circ\pi\circ(\alpha_i|_V,\beta)=r_Y\circ\beta$, irrespective of $i\in\lbrace 1,2\rbrace$. We also note that the codomain of $\Omega^\alpha|_V$ is exactly the \emph{do}main of the inverse $\Psi=(\Psi_1,\Psi_2)$ of the action map of the balanced tensor product $G\LAalong{L_X}X\tensor_H Y\xrightarrow{R_Y}H_0$ (given explicitly in \cref{proposition:pre-principal balanced tensor product}). In particular we then get a smooth map
		\begin{equation*}
			\Psi_1\circ\Omega^\alpha|_V:
			V
			\xrightarrow{\quad\Omega^\alpha|_V\quad}
			\left(X\tensor_H Y\right)\times_{G_0}^{R_Y,R_Y}\left(X\tensor_H Y\right)
			\xrightarrow{\quad\Psi_1\quad}
			G.
		\end{equation*}
		We now extend this map to the entire domain $U_\alpha$, and show that it is independent on the choice of plot $\beta$. For that, pick two points $t,\overline{t}\in U_\alpha$, so that by subductiveness of the left moment map $l_Y$ we can find two plots, $\beta:V\to Y$ and $\overline{\beta}:\overline{V}\to Y$, defined on open neighbourhoods of $t$ and $\overline{t}$, respectively, such that $r_X\circ \alpha_i|_V=l_Y\circ\beta$ and $r_X\circ \alpha_i|_{\overline{V}}=l_Y\circ\overline{\beta}$. Following the above construction, we get two smooth maps:
		\begin{align*}
			\Omega^\alpha|_V
			:s
			&\longmapsto 
			\left(\alpha_1|_V(s)\tensor \beta(s),\alpha_2|_V(s)\tensor\beta(s)\right),\\
			\overline{\Omega}^\alpha|_{\overline{V}}:
			s
			&\longmapsto
			\left(\alpha_1|_{\overline{V}}(s)\tensor\overline{\beta}(s),
			\alpha_2|_{\overline{V}}(s)\tensor\overline{\beta}(s)\right).	
		\end{align*}
		We now remark an important characterisation of $\Psi$, as a consequence of it being a diffeomorphism and inverse to the action map. Namely, $\Psi_1(x_1\tensor y_1,x_2\tensor y_2)\in G$ is the \emph{unique} arrow $g\in G$ satisfying $gx_2\tensor y_2=x_1\tensor y_1$. Therefore, $\Psi_1\circ\Omega^\alpha|_V(s)\in G$ is the unique arrow such that
		\begin{equation*}
			\left[\Psi_1\circ\Omega^\alpha|_V(s)\right]\cdot\left(\alpha_2|_V(s)\tensor\beta(s)\right)= \alpha_1|_V(s)\tensor\beta(s).
		\end{equation*}
		By \cref{lemma:actions of weakly invertible bibundle are free} all of the four actions of the original bibundles are free. Consequently, applying \cref{lemma:free action and balanced tensor product}, since the second component in each term is just $\beta(s)$, this means that $\Psi_1\circ\Omega^\alpha|_V(s)$ is the unique arrow in $G$ such that
		\begin{equation*}
			\Psi_1\circ\Omega^\alpha|_V(s)\cdot\alpha_2|_V(s)=\alpha_1|_V(s),
		\end{equation*}
		where the tensor with $\beta(s)$ can be removed. But, for exactly the same reasons, if we take $s\in V\cap \overline{V}$, then $\Psi_1\circ\overline{\Omega}^\alpha|_{\overline{V}}(s)\in G$ is \emph{also} the unique arrow such that
		\begin{equation*}
			\Psi_1\circ\overline{\Omega}^\alpha|_{V\cap \overline{V}}(s)\cdot\alpha_2|_{V\cap\overline{V}}(s)=\alpha_1|_{V\cap\overline{V}}(s),
		\end{equation*}
		proving that 
		\begin{equation*}
			\Psi_1\circ\Omega^\alpha|_{V\cap\overline{V}}=\Psi_1\circ\overline{\Omega}^\alpha|_{V\cap\overline{V}}.
		\end{equation*}
		This shows that on the overlaps $V\cap \overline{V}$ the map $\Psi_1\circ\Omega^\alpha|_{V\cap\overline{V}}$ does \emph{not} depend on the plots $\beta$ and $\overline{\beta}$. This allows us to extend $\Psi_1\circ\Omega^\alpha|_V$, in a well-defined way, to the entire domain of $U_\alpha$. We do this as follows. For every $t\in U_\alpha$ there exists a plot $\beta_t:V_t\to Y$, defined on an open neighbourhood $V_t\ni t$, such that $r_X\circ\alpha_i|_{V_t}=l_Y\circ\beta_t$. Clearly, this gives an open cover $(V_t)_{t\in U_\alpha}$ of $U_\alpha$. For $t\in U_\alpha$ we then set $\Psi_1\circ\Omega^\alpha(t):=\Psi_1\circ\Omega^\alpha|_{V_t}(t)$. Hence we get a well-defined function $\Psi_1\circ\Omega^\alpha:U_\alpha\to G$, which is smooth by the Axiom of Locality.		
	\end{construction}
	
	The main observation now is that, as the plot $\alpha$ is centred at $(x_1,x_2)$, we get that $\Psi_1\circ\Omega^\alpha(0)$ is the unique arrow in $G$ such that $\Psi_1\circ\Omega^\alpha(0)\cdot x_2=x_1$. This is exactly the property that characterises the division $\langle x_1,x_2\rangle_G$!
	
	\begin{proposition}
		\label{proposition:weakly invertible is pre-biprincipal}
		A weakly invertible diffeological bibundle is pre-biprincipal.
		
		\begin{proof}
			The bulk of the work has been done in \cref{construction:local division map}. Start with a diffeological bibundle $G\LAalong{l_X}X\RAalong{r_X}H$ and a weak inverse $H\LAalong{l_Y}Y\RAalong{r_Y}G$. We shall define a smooth division map $\langle\cdot,\cdot\rangle_G$ for the left $G$-action. For $(x_1,x_2)\in X\times_{H_0}^{r_X,r_X}X$, we know by the Axiom of Covering that the constant map $\const_{(x_1,x_2)}:\mathbb{R}\to X\times_{H_0}^{r_X,r_X}X$ is a plot centred at $(x_1,x_2)$. We use the shorthand $\Psi_1\circ \Omega^{(x_1,x_2)}$ to denote the map $\Psi_1\circ\Omega^\alpha$ defined by the plot $\alpha=\const_{(x_1,x_2)}$, and then write:
			\begin{equation*}
				\langle x_1,x_2\rangle_G:= \Psi_1\circ\Omega^{(x_1,x_2)}(0).
			\end{equation*}
			That just leaves us to show that this map is smooth. For that, take an arbitrary plot $\alpha:U_\alpha\to X\times_{H_0}^{r_X,r_X}X$ of the fibred product. We need to show that $\langle \cdot,\cdot\rangle_G\circ \alpha$ is a plot of $G$. For any $t\in U_\alpha$, we have that
			\begin{equation*}
				\langle \alpha_1(t),\alpha_2(t)\rangle_G=\Psi_1\circ \Omega^{\alpha(t)}(0)
			\end{equation*}
			is the unique arrow in $G$ such that 
			\begin{equation*}
				\Psi_1\circ\Omega^{\alpha(t)}(0)\cdot\const^2_{\alpha(t)}(0)=\const^1_{\alpha(t)}(0),
			\end{equation*}
			where $\const^i$ denotes the $i$th component of the constant plot. But then $\const^i_{\alpha(t)}(0)=\alpha_i(t)$, and we already know that $\Psi_1\circ\Omega^\alpha(t)\in G$ is the unique arrow that sends $\alpha_2(t)$ to $\alpha_1(t)$, so we have:
			\begin{equation*}
				\Psi_1\circ\Omega^{\alpha(t)}(0)=\Psi_1\circ\Omega^\alpha(t),\qquad\text{which means}\qquad \langle\cdot,\cdot\rangle_G\circ\alpha=\Psi_1\circ\Omega^\alpha.
			\end{equation*}
			But the right hand side $\Psi_1\circ\Omega^\alpha:U_\alpha\to G$ is a plot of $G$ as per \cref{construction:local division map}, proving that the map $\langle\cdot,\cdot\rangle_G$ is smooth. It is quite evident from its construction that it satisfies exactly the properties of a division map, and it is now easy to verify that
			\begin{equation*}
				\left(\langle\cdot,\cdot\rangle_G,\pr_2|_{X\times_{H_0}X}\right):X\times_{H_0}^{r_X,r_X}X\longrightarrow G\times_{G_0}^{\source,l_X}X				
			\end{equation*}
			is a smooth inverse of the action map (see \cref{section:division map}). The fact that it lands in the right codomain, i.e., $\source(\langle x_1,x_2\rangle_G)=l_X(x_2)$, follows from the properties of $\Psi$ as the inverse of the action map of the balanced tensor product. Therefore $G\LAalong{l_X}X\xrightarrow{r_X}H_0$ is a pre-principal bundle. An analogous argument will show that $G_0\xleftarrow{l_Y}X\RAalong{r_X}H$ is also pre-principal, and hence we have proved the claim.
		\end{proof}
	\end{proposition}
	
	We can now prove our main theorem:
	\begin{theorem}
		\label{theorem:weakly invertible bibundles are the biprincipal ones}
		A bibundle is weakly invertible in $\DiffBiBund$ if and only if it is biprincipal. That means: two diffeological groupoids are Morita equivalent if and only if they are equivalent in $\DiffBiBund$.
		
		\begin{proof}
			One of the implications is just \cref{proposition:biprincipal bibundle is weakly invertible}. The other now follows from a combination of \cref{proposition:weakly invertible implies bisubductive,proposition:weakly invertible is pre-biprincipal}.
		\end{proof}
	\end{theorem}
	
	This significantly generalises \cite[Proposition 4.21]{landsman2001quantized}, not only in that we have a generalisation to a diffeological setting, but also in that it considers a more general type of bibundle. It justifies the bicategory $\DiffBiBund$ as being the appropriate setting for Morita equivalence of diffeological groupoids. It also shows that the assumptions of left principality of the Lie groupoid bibundles appear to be more like technical necessities for getting a well defined bicategory of Lie groupoids and bibundles, rather than being meaningful assumptions on the underlying smooth structure of the bibundles. In \cref{section:diffeological bibundles between Lie groupoids} we discuss other aspects of diffeological Morita equivalence between Lie groupoids. A possible \emph{category of fractions} approach to Morita equivalence of diffeological groupoids is discussed in \cite[Chapter V]{schaaf2020diffeology-groupoids-and-ME}.

\section{Some Morita Invariants}
	\label{section:some applications}
	In theories of Morita equivalence, there are often interesting properties that are naturally Morita invariant. In this section we discuss some results that generalise several well known Morita invariants of Lie groupoids to the diffeological setting. These include: invariance of the orbit spaces (\cref{definition:groupoid orbit space}), of being \emph{fibrating} (\cref{definition:fibrating groupoid}), and of the action categories (\cref{definition:action category})  The proofs are taken from \cite[Chapter IV]{schaaf2020diffeology-groupoids-and-ME}.
	
	\subsection{Invariance of orbit spaces}
	\label{section:invariance orbit spaces}
	It is a well known result that if two Lie groupoids $\grpd{G}$ and $\grpd{H}$ are Morita equivalent (in the Lie groupoid sense), then there is a \emph{homeo}morphism between their orbit spaces $G_0/G$ and $H_0/H$ \cite[Lemma 2.19]{crainic2018orbispaces}. The following theorem shows that, for diffeological groupoids, we get a genuine \emph{diffeo}morphism. The construction of the underlying function is the same as for the Lie groupoid case, which is sketched in the proof of \cite[Lemma 2.19]{crainic2018orbispaces}, and which we describe below in detail.
	
	\begin{theorem}
		\label{theorem:morita equivalent groupoids have diffeomorphic orbit spaces (bibundle proof)}
		If $\grpd{G}$ and $\grpd{H}$ are two Morita equivalent diffeological groupoids, then there is a diffeomorphism $G_0/G\cong H_0/H$ between their orbit spaces.
		
		\begin{proof}
			Let $G\LAalong{l_X}X\RAalong{r_X}H$ be the bibundle instantiating the Morita equivalence. Our first task will be to construct a function $\Phi: G_0/G\to H_0/H$ between the orbit spaces. The idea is to lift a point $a\in G_0$ of the base of the groupoid to its $l_X$-fibre, which by right principality is just an $H$-orbit in $X$, and then to project this orbit down to the other base $H_0$ along the right moment map $r_X$. The fact that the bundle is biprincipal ensures that this can be done in a consistent fashion.
			
			We are dealing with \emph{four} actions here, so we need to slightly modify our notation to avoid confusion. If $a\in G_0$ is an object in the groupoid $G$, we shall denote its orbit by $\Orb_{G_0}(a)$, which, as usual, is just the set of all points $a'\in G_0$ such that there exists an arrow $g:a\to a'$ in $G$. Similarly, for $b\in H_0$ we write $\Orb_{H_0}(b)$. On the other hand, we have two actions on $X$, for whose orbits we use the standard notations $\Orb_G(x)$ and $\Orb_H(x)$, where $x\in X$.
			
			Now, start with a point $a\in G_0$, and consider its fibre $l_X^{-1}(a)$ in $X$. Since the bibundle is right subductive, the map $l_X$ is surjective, so this fibre is non-empty and we can find a point $x_a\in l_X^{-1}(a)$. We claim that the expression $\Orb_{H_0}\circ r_X(x_a)$ is independent on the choice of the point $x_a$ in the fibre. For that, take another point $x_a'\in l_X^{-1}(a)$. This gives the equation $l_X(x_a)=l_X(x_a')$, and since bibundle is right pre-principal, we get a unique arrow $h\in H$ such that $x_a'=x_ah$. From the definition of a right groupoid action, this in turn gives the equations $r_X(x_a')= \source(h)$ and $r_X(x_a)=\target(h)$, which proves the claim. To summarise, whenever $x_a,x_a'\in l_X^{-1}(a)$ are two points in the same $l_X$-fibre, then we have:
			\begin{equation}
			\label{equation:independent H_0 orbit on point in fibre}
			\Orb_{H_0}\circ r_X(x_a)
			=
			\Orb_{H_0}\circ r_X(x_a').
			\tag{{\color{RadboudRed}$\clubsuit$}}
			\end{equation}
			Next we want to show that neither is this expression dependent on the point $a\in G_0$, but rather on its orbit $\Orb_{G_0}(a)$. For this, take another point $b\in \Orb_{G_0}(a)$, so there exists some arrow $g:a\to b$ in $G$. Pick then $x\in l_X^{-1}(a)$ and $y\in l_X^{-1}(b)$. This means that $\source(g)=l_X(x)$ and $\target(g)=l_X(y)$, which means that if we let $g$ act on the point $x$ we get a point $gx\in l_X^{-1}(b)$, in the same $l_X$-fibre as $y$. Then using equation \eqref{equation:independent H_0 orbit on point in fibre} applied to $gx$ and $y$, and the $G$-invariance of the right moment map $r_X$, we immediately get:
			\begin{equation*}
				\Orb_{H_0}\circ r_X(x) = \Orb_{H_0}\circ r_X(gx) = \Orb_{H_0}\circ r_X(y).
			\end{equation*}
			Using this, we can now conclude that there is a well-defined function
			\begin{equation*}
				\Phi
				:
				G_0/G
				\longrightarrow
				H_0/H;
				\qquad
				\Orb_{G_0}(a)
				\longmapsto
				\Orb_{H_0}\circ r_X(x_a),
			\end{equation*}
			that is neither dependent on the point $a$ in the orbit $\Orb_{G_0}(a)$, nor on the choice of the point $x_a\in l_X^{-1}(a)$ in the fibre. Note that this function exists by virtue of right subductivity (and the Axiom of Choice), which ensures that the left moment map $l_X$ is a surjection (and for each $a$ there exists an $x_a$).
			
			Either by replacing $G\LAalong{l_X}X\RAalong{r_X}H$ by its opposite bibundle, or by switching the words `left' and `right', the above argument analogously gives a function going the other way:
			\begin{equation*}
				\Psi
				:
				H_0/H
				\longrightarrow
				G_0/G;
				\qquad
				\Orb_{H_0}(b)
				\longmapsto
				\Orb_{G_0}\circ l_X(y_b),
			\end{equation*}
			where now $y_b\in r_X^{-1}(b)$ is some point in the fibre of the right moment map $r_X$. We claim that $\Phi$ and $\Psi$ are mutual inverses. To see this, pick a point $a\in G_0$, a point $x_a\in l_X^{-1}(a)$, a point $y_{r_X(x_a)}\in r_X^{-1}(r_X(x_a))$. Then we can write
			\begin{equation*}
				\Psi\circ \Phi\left(
				\Orb_{G_0}(a)
				\right)
				=
				\Psi\left(
				\Orb_{H_0}(r_X(x_a))
				\right)
				=
				\Orb_{G_0}\left(l_X(y_{r_X(x_a)})\right).
			\end{equation*}
			We also have, by choice, the equation $r_X(x_a)=r_X(y_{r_X(x_a)})$, so by left pre-principality there exists an arrow $g\in G$ such that $gx_a=y_{r_X(x_a)}$. By definition of a left groupoid action, this then further gives 
			\begin{equation*}
				\source(g)=l_X(x_a)=a
				\qquad\text{and}\qquad
				\target(g)=l_X(y_{r_X(x_a)}).
			\end{equation*}
			This proves that the right-hand side of the previous equation is equal to
			\begin{equation*}
				\Orb_{G_0}\left(l_X(y_{r_X(x_a)})\right)
				=
				\Orb_{G_0}(a),
			\end{equation*}
			which gives $\Psi\circ\Phi=\id_{G_0/G}$. Through a similar argument, using right pre-principality, we obtain that $\Phi\circ\Psi=\id_{H_0/H}$. 
			
			To finish the proof, it suffices to prove that both $\Phi$ and $\Psi$ are smooth. Again, due to the symmetry of the situation, and since the bibundle $G\LAalong{l_X}X\RAalong{r_X}H$ is biprincipal, we shall only prove that $\Phi$ is smooth. The proof for $\Psi$ will follow analogously. Since $\Orb_{G_0}$ is a subduction, to prove that $\Phi$ is smooth it suffices by \cref{lemma:properties of subductions}\emph{(3)} to prove that $\Phi\circ\Orb_{G_0}$ is smooth. Since the left moment map $l_X$ is a surjection, using the Axiom of Choice we pick a section $\sigma:G_0\to X$, which replaces our earlier notation of $\sigma(a)=:x_a$. From the way $\Phi$ is defined, we see that we get a commutative diagram:
			\begin{equation*}
				\begin{tikzcd}
					G_0 \arrow[d, "\Orb_{G_0}"'] \arrow[r, "\sigma"] & X \arrow[r, "r_X"] & H_0 \arrow[d, "\Orb_{H_0}"] \\
					G_0/G \arrow[rr, "\Phi"']                        &                    & H_0/H      .               
				\end{tikzcd}
			\end{equation*}
			We are therefore to show that $\Orb_{H_0}\circ r_X\circ \sigma$ is smooth. For this, pick a plot $\alpha:U_\alpha\to G_0$ of the base space. By right subductivity, the left moment map $l_X$ is a subduction, so locally $\alpha|_V=l_X\circ \beta$, where $\beta$ is some plot of $X$. Now, note that, for all $t\in V$, both the points $\beta(t)$ and $\sigma\circ l_X\circ \beta(t)$ are elements of the fibre $l_X^{-1}(l_X\circ\beta(t))$. Therefore, by equation \eqref{equation:independent H_0 orbit on point in fibre} we get:
			\begin{equation*}
				\Orb_{H_0}\circ r_X\circ \sigma\circ \alpha|_V
				=
				\Orb_{H_0}\circ r_X\circ \sigma\circ l_X\circ \beta
				=
				\Orb_{H_0}\circ r_X\circ \beta.
			\end{equation*}
			The right-hand side of this equation is clearly smooth (and no longer dependent on the choice of section $\sigma$). By the Axiom of Locality for $G_0$, it follows that $\Orb_{H_0}\circ r_X\circ \sigma\circ \alpha$ is globally smooth, and since the plot $\alpha$ was arbitrary, this proves that $\Phi\circ\Orb_{G_0}$ is smooth. Hence, $\Phi$ is smooth. After an analogous argument that shows $\Psi$ is smooth, the desired diffeomorphism between the orbit spaces follows.
		\end{proof}
	\end{theorem}

	\subsection{Invariance of fibration}
	\label{section:invariance of fibration}
	The theory of diffeological (principal) fibre bundles is shown in \cite[Chapter 8]{iglesias2013diffeology} to be fully captured by the following notion:
	\begin{definition}
		\label{definition:fibrating groupoid}
		A diffeological groupoid $\grpd{G}$ is called \emph{fibrating} (or a \emph{fibration groupoid}) if the \emph{characteristic map} $(\target,\source):G\to G_0\times G_0$ is a subduction.
	\end{definition}
	
	This leads to a theory of diffeological fibre bundles that is able to treat the standard smooth locally trivial (principal) fibre bundles of smooth manifolds, but also bundles that are not (and could not meaningfully be) locally trivial. It is then natural to ask if this property of diffeological groupoids is invariant under Morita equivalence. The following theorem proves that this is the case:
	
	\begin{theorem}
		Let $\grpd{G}$ and $\grpd{H}$ be two Morita equivalent diffeological groupoids. Then $\grpd{G}$ is fibrating if and only if $\grpd{H}$ is fibrating.
		
		\begin{proof}
			Because Morita equivalence is an equivalence relation, it suffices to prove that if $\grpd{G}$ is fibrating, then so is $\grpd{H}$. Denoting the characteristic maps of these groupoids by $\chi_G=\left(\target_G,\source_G\right)$ and $\chi_H=\left(\target_H,\source_H\right)$, assume that $G$ is fibrating, so that $\chi_G$ is a subduction. Our goal is to show $\chi_H$ is also a subduction.
			
			To begin with, take an arbitrary plot $\alpha=(\alpha_1,\alpha_2):U_\alpha \to H_0\times H_0$, and fix an element $t\in U_\alpha$. We thus need to find a plot $\Phi:W\to H$, defined on an open neighbourhood $t\in W\subseteq U_\alpha$, such that $\alpha|_W=\chi_H\circ\Phi$. Morita equivalence yields a biprincipal bibundle $G\LAalong{l_X}X\RAalong{r_X}H$. To construct the plot $\Phi$, we use almost all of the structure of this bibundle. 
			
			The right moment map $r_X:X\to H_0$ is a subduction, so for each of the components $\alpha_i$ of $\alpha$ we get a plot $\beta_i:U_i\to X$, defined on an open neighbourhood $t\in U_i\subseteq U_\alpha$, such that $\alpha_i|_{U_i}=r_X\circ\beta_i$. Define $U:=U_1\cap U_2$, which is another open neighbourhood of $t\in U_\alpha$, and introduce the notation
			\begin{equation*}
				\beta:=(\beta_1|_U,\beta_2|_U):U\longrightarrow X\times X.
			\end{equation*}
			Composing with the left moment map $l_X:X\to G_0$, we get $(l_X\times l_X)\circ\beta:U\to G_0\times G_0$. It is here that we use that $\grpd{G}$ is fibrating. Because of that, we can find an open neighbourhood $t\in V\subseteq U\subseteq U_\alpha$ and a plot $\Omega:V\to G$ such that 
			\begin{equation}
			\label{equation:chi_G Omega}
			\chi_G\circ\Omega
			=
			\left.(l_X\times l_X)\circ\beta\right|_V.
			\tag{{{\color{RadboudRed}$\spadesuit$}}}
			\end{equation}
			This means that $\target_G\circ\Omega = l_X\circ\beta_1|_V$ and $\source_G \circ\Omega = l_X\circ \beta_2|_V$. Let ${\varphi_G:X\tensor_H \overline{X}\to G}$ be the biequivariant diffeomorphism from \cref{proposition:biprincipal bibundle is weakly invertible}. Using the plot $\Omega$ we just obtained, we get another plot $\varphi_G^{-1}\circ\Omega:V\to X\tensor_H\overline{X}$. Now, since the canonical projection $\pi_H:X\times_{H_0}^{r_X,r_X}\overline{X}\to X\tensor_H\overline{X}$ of the diagonal $H$-action is a subduction, we can find an open neighbourhood $t\in W\subseteq V$ and a plot $\omega:W\to X\times_{H_0}^{r_X,r_X}\overline{X}$ such that
			\begin{equation}
			\label{equation:pi_H omega}
			\pi_H\circ \omega = \varphi_G^{-1}\circ\Omega|_W.
			\tag{{\color{RadboudRed}$\clubsuit$}}
			\end{equation}
			Note that the plot $\omega$ decomposes into its components $\omega_1,\omega_2:W\to X$, which satisfy $r_X\circ\omega_1=r_X\circ\omega_2$. Using the biequivariance of $\varphi_G$ and the defining relation $L_X\circ\pi_H = l_X\circ \pr_1|_{X\times_{H_0}\overline{X}}$ we find:
			\begin{equation*}
				l_X\circ\beta_1|_W
				=
				\target_G\circ\Omega|_W
				=
				L_X\circ\varphi_G^{-1}\circ\Omega|_W
				=
				L_X\circ \pi_H\circ\omega
				=
				l_X\circ\pr_1|_{X\times_{H_0}\overline{X}}\circ\omega
				=
				l_X\circ\omega_1,
			\end{equation*}
			where the first equality follows from the equation \eqref{equation:chi_G Omega}, and the third one from \eqref{equation:pi_H omega}. Similarly, we find $l_X\circ\beta_2|_W = l_X\circ\omega_2$. These two equalities give two well-defined plots, one for each $i\in\lbrace 1,2\rbrace$, given by
			\begin{equation*}
				\beta_i|_W\tensor\omega_i:=\pi_G\circ\left(\beta_i|_W,\omega_i\right)
				:W
				\xrightarrow{\quad (\beta_i|_W,\omega_i)\quad} 
				\overline{X}\times_{G_0}^{l_X,l_X}X
				\xrightarrow{\quad\pi_G\quad}
				\overline{X}\tensor_G X,
			\end{equation*}
			where $\pi_G:\overline{X}\times_{G_0}^{l_X,l_X}X \to \overline{X}\tensor_G X$ is the canonical projection of the diagonal $G$-action. We can now apply the biequivariant diffeomorphism $\varphi_H:\overline{X}\tensor_G X\To H$ from \cref{proposition:biprincipal bibundle is weakly invertible} to get two plots in $H$. It is from these two plots that we will create $\Phi$. Here it is absolutely essential that we have constructed the plot $\omega$ such that $r_X\circ\omega_1=r_X\circ\omega_2$, because that means that the sources of these two plots in $H$ will be equal, and hence they can be composed if we first invert one of them component-wise. To see this, use the biequivariance of $\varphi_H$ to calculate
			\begin{equation*}
				\source_H\circ\varphi_H\circ\left(\beta_i|_W\tensor\omega_i\right)
				=
				R_X\circ\left(\beta_i|_W\tensor\omega_i\right)
				=
				r_X\circ\pr_2|_{\overline{X}\times_{G_0}X}\circ\left(\beta_i|_W,\omega_i\right)
				=r_X\circ\omega_i,
			\end{equation*}
			and similarly:
			\begin{equation*}
				\target_H\circ\varphi_H\circ\left(\beta_i|_W\tensor\omega_i\right)
				=
				L_{\overline{X}}\circ\left(\beta_i|_W\tensor\omega_i\right)
				=
				r_X\circ\pr_1|_{\overline{X}\times_{G_0}X}\circ\left(\beta_i|_W,\omega_i\right)
				=
				r_X\circ \beta_i|_W
				= \alpha_i|_W.
			\end{equation*}
			Of course, if we switch $\beta_i|_W\tensor\omega_i$ to $\omega_i\tensor \beta_i|_W$, which is defined in the obvious way, then the right-hand sides of the above two equations will switch. So, for every $s\in W$, the expression $\varphi_H\left(\omega_2(s)\tensor \beta_2(s)\right)$ is an arrow in $H$ from $r_X\circ\beta_2(s)=\alpha_2(s)$ to $r_X\circ\omega_2(s)$, and $\varphi_H\left(\beta_1(s)\tensor\omega_1(s)\right)$ is an arrow from $r_X\circ\omega_1(s)=r_X\omega_2(s)$ to $r_X\circ\beta_1(s)=\alpha_1(s)$, which can hence be composed to give an arrow from $\alpha_2(s)$ to $\alpha_1(s)$. This is exactly the kind of arrow we want. Therefore, for every $s\in W$, we get a commutative triangle in the groupoid $H$, which defines for us the plot $\Phi:W\to H$:
			\begin{equation*}
				\begin{tikzcd}[column sep = tiny]
					\alpha_2(s) \arrow[rd, "\varphi_H\left(\omega_2(s)\tensor\beta_2(s) \right)"'] \arrow[rr, "\Phi(s)", dashed] &                                                                                    & \alpha_1(s) \\
					& r_X\circ\omega_1(s). \arrow[ru, "\varphi_H\left(\beta_1(s)\tensor\omega_1(s)\right)"']
				\end{tikzcd}
			\end{equation*}
			The map $\Phi$ is clearly smooth, because inversion and multiplication in $H$ are smooth. Hence we have defined the plot $\Phi$, and by the above diagram it is clear that it satisfies
			\begin{equation*}
				\chi_H\circ\Phi = (\target_H\circ\Phi,\source_H\circ\Phi)=\alpha|_W.
			\end{equation*}
			Thus we may at last conclude that $\chi_H$ is a subduction, and hence that $\grpd{H}$ is also fibrating.
		\end{proof}
	\end{theorem}
	
	\subsection{Invariance of representations}
	\label{section:invariance of representations}
	In the Morita theory of rings, it holds that two rings are Morita equivalent if and only if their categories of modules are equivalent. For groupoids, even discrete ones, this is no longer an ``if and only if'' proposition, but merely an ``only if''. Nevertheless, it is known that the result transfers to Lie groupoids as well \cite[Theorem 6.6]{landsman2001bicategories}, and here we shall prove that it transfers also to diffeology.
	
	\begin{theorem}
		\label{theorem:morita equivalent groupoids have equivalent action categories}
		Suppose that $\grpd{G}$ and $\grpd{H}$ are Morita equivalent diffeological groupoids. Then the action categories $\Act(\grpd{G})$ and $\Act(\grpd{H})$ are categorically equivalent.
		
		\begin{proof}
			If $\grpd{G}$ and $\grpd{H}$ are Morita equivalent, there exists a biprincipal bibundle $G\LAalong{l_X}X\RAalong{r_X} H$. Recall from \cref{definition:action category} the notion of action categories and from \cref{definition:induced action functor} that of induced action functors. We claim that 
			\begin{align*}
				X\tensor_H - :\Act(\grpd{H})&\longrightarrow\Act(\grpd{G}),
				\\
				\overline{X}\tensor_G-:\Act(\grpd{G})&\longrightarrow \Act(\grpd{H})
			\end{align*}
			are mutually inverse functors up to natural isomorphism. To see this, take a left $H$ action $H\LAalong{l_Y}Y$. Then 
			\begin{equation*}
				\left(\overline{X}\tensor_G-\right)\circ\left(X\tensor_H-\right)[H\LAalong{l_Y} Y]
				=
				\left(\overline{X}\tensor_G -\right)\left[G\LAalong{L_X} X\tensor_HY\right]
				=
				H\LAalong{L_{\overline{X}}}\left(\overline{X}\tensor_G\left(X\tensor_HY\right)\right).
			\end{equation*}
			Therefore, we need to construct a natural biequivariant diffeomorphism
			\begin{equation*}
				\mu_Y:\overline{X}\tensor_G\left(X\tensor_H Y\right)\longrightarrow Y.
			\end{equation*}
			For this, we collect the biequivariant diffeomorphisms from \cref{proposition:identity bibundle is weak identity,proposition:associativity of balanced tensor product,proposition:biprincipal bibundle is weakly invertible}. Let us denote them by
			\begin{align*}
				&A_Y:\overline{X}\tensor_G\left(X\tensor_H Y\right)
				\longrightarrow 
				\left(\overline{X}\tensor_G X\right)\tensor_HY,
				\\
				&\varphi_H:\overline{X}\tensor_G X\longrightarrow H,
				\\
				&M_{Y}:H\tensor_H Y\longrightarrow Y,
			\end{align*}
			describing the association up to isomorphism, the division map of the bibundle, and the left action $H\laction Y$, respectively. We then define
			\begin{equation*}
				\mu_Y:= M_Y\circ\left(\varphi_H\tensor\id_Y\right)\circ A_Y.
			\end{equation*}
			Note that $(\varphi_H\tensor\id_Y)$ is still a biequivariant diffeomorphism. The naturality square of the natural transformation ${\mu:\left(\overline{X}\tensor_G-\right)\circ\left(X\tensor_H-\right)\To \id_{\Act(H)}}$ then becomes:
			\begin{equation*}
				\begin{tikzcd}[column sep =large]
					\overline{X}\tensor_G\left(X\tensor_H Y\right) \arrow[r, "\mu_Y"] \arrow[d, "\id_{\overline{X}}\tensor(\id_X\circ\varphi)"'] & Y \arrow[d, "\varphi"] \\
					\overline{X}\tensor_G\left(X\tensor_H Z\right) \arrow[r, "\mu_Z"']                                                           & {Z,}                  
				\end{tikzcd}
			\end{equation*}
			where $\varphi:Y\to Z$ is an $H$-equivariant smooth map. It follows from the structure of these maps that the naturality square commutes. The top right corner of the diagram becomes:
			\begin{align*}
				\varphi\circ\mu_Y\left(x_1\tensor(x_2\tensor y)\right)
				&=
				\varphi\circ M_Y\circ(\varphi_H\tensor\id_Y)\circ A_Y\left(x_1\tensor(x_2\tensor y)\right)
				\\&=
				\varphi\circ M_Y\circ(\varphi_H\tensor\id_Y)\left((x_1\tensor x_2)\tensor y\right)
				\\&=
				\varphi\circ M_Y\left(\varphi_H(x_1\tensor x_2)\tensor y\right)
				\\&=
				\varphi\left(\varphi_H(x_1\tensor x_2)y\right)
				\\&=
				\varphi_H(x_1\tensor x_2)\varphi(y),
			\end{align*}
			where the very last step follows from $H$-equivariance of $\varphi$. Following a similar calculation, the bottom left corner evaluates as
			\begin{align*}
				\mu_Z\circ\left(\id_{\overline{X}}\tensor(\id_X\tensor\varphi)\right)
				&=
				M_Z\circ(\varphi_H\tensor\id_Z)\circ A_Z\circ\left(\id_{\overline{X}}\tensor(\id_X\tensor\varphi)\right)
				\\&=
				M_Z\circ(\varphi_H\tensor\id_Z)\circ\left((\id_{\overline{X}}\tensor\id_X)\tensor\varphi\right)
				\\&=
				M_Z\circ(\varphi_H\tensor\varphi),
			\end{align*}
			which, when evaluated, gives exactly the same as the above expression for the top right corner. This proves that $\mu$ is natural, and since every of its components is an $H$-equivariant diffeomorphism, it follows that $\mu$ is a natural isomorphism. The fact that the composition $\left(X\tensor_H -\right)\circ\left(\overline{X}\tensor_G-\right)$ is naturally isomorphic to $\id_{\Act(G)}$ follows from an analogous argument. Hence the categories $\Act(\grpd{G})$ and $\Act(\grpd{H})$ are equivalent, as was to be shown.
		\end{proof}
	\end{theorem}

\section{Discussion and Suggestions for Future Research}
	\label{section:closing section}
	\subsection{Diffeological bibundles between Lie groupoids}
	\label{section:diffeological bibundles between Lie groupoids}
	As we saw in \cref{example:Lie ME is also diffeological ME}, if two Lie groupoids are \emph{Lie} Morita equivalent (i.e. Morita equivalent in the Lie groupoid sense \cite[Definition 2.15]{crainic2018orbispaces}), then they are also \emph{diffeologically} Morita equivalent. This is simply due to the fact that surjective submersions between smooth manifolds are in particular also subductions, and hence a Lie principal groupoid bundle is also diffeologically principal. But, what if $\grpd{G}$ and $\grpd{H}$ are two \emph{Lie} groupoids, such that there exists a \emph{diffeological} biprincipal bibundle $G\LAalong{l_X}X\RAalong{r_X}H$ between them. What does that say about the \emph{Lie} Morita equivalence of $G$ and $H$? This still remains an open question (\cref{question:does inclusion pseudofunctor reflect weak equivalence}). In this section we discuss some related results, which also pertain to our choice of subductions over \emph{local} subductions for the development of the general theory. A slightly more detailed discussion is in \cite[Section 4.4.3]{schaaf2020diffeology-groupoids-and-ME}. In light of \cref{proposition:local subductions are surjective submersions}, the source and target maps of a Lie groupoid are local subductions (cf. \cref{proposition:source map is subduction}), and we can therefore introduce the following class of diffeological groupoids:
	
	\begin{definition}
		\label{definition:locally subductive groupoids}
		We say a diffeological groupoid $\grpd{G}$ is \demph{locally subductive} if its source and target maps are local subductions\footnote{It would be tempting to call such groupoids \emph{``diffeological Lie groupoids,''} but this would conflict with earlier established terminology of so-called \emph{diffeological Lie groups} in \cite[Article 7.1]{iglesias2013diffeology} and \cite{leslie2003diffeological,magnot2018group}.}. Clearly, every Lie groupoid is a locally subductive diffeological groupoid.
	\end{definition}
	Looking at the structure of the proofs in \cref{section:diffeological groupoid actions and bundles,section:diffeological bibundles}, it appears as if they can be generalised to a setting where we replace all subductions by local subductions. In doing so, we would get a theory of locally subductive groupoids, locally subductive groupoid bundles, and the corresponding notions for bibundles and Morita equivalence, which, as it appears, would follow the same story as we have so far presented. An upside to that framework would be that it directly returns the original theory of Morita equivalence for Lie groupoids, once we restrict our diffeological spaces to smooth manifolds. In this section we shall prove that, even in the slightly more general setting of \cref{section:diffeological bibundles}, the diffeological bibundle theory reduces to the Lie groupoid theory in the correct way. We do this by proving that the moment maps of a biprincipal bibundle between locally subductive groupoids have to be local subductions as well (\cref{lemma:biprincipal bib between loc subd groupoids is locally bisubductive}). In hindsight, this provides more justification for our choice of starting with subductions instead of local subductions.  One consequence of this choice is that it allows for groupoid bundles that are truly \emph{pseudo}-bundles, in the sense of \cite{pervova2016diffeological}. The notion of pseudo-bundles seems to be the correct notion in the setting of diffeology to generalise all bundle constructions on manifolds, at least if we want to treat (internal) tangent bundles as such (see \cite{christensen2016tangent}). There exists diffeological spaces whose internal tangent bundle is not a local subduction \cite[Example 3.17]{christensen2016tangent}. If we had defined principality of a groupoid bundle to include \emph{local} subductiveness, these examples would not be treatable by our theory of Morita equivalence.
	
	\begin{lemma}
		\label{lemma:projection onto balanced tensor product is local subduction}
		Let $G\LAalong{l_X}X\RAalong{r_X}H$ be a diffeological bibundle, where $\grpd{H}$ is a locally subductive groupoid. Then the canonical projection map $\pi_H:X\times_{H_0}^{r_X,r_X}\overline{X}\to X\tensor_H\overline{X}$ is a local subduction.
		
		\begin{proof}
			Let $\alpha:(U_\alpha,0)\to (X\tensor_H\overline{X},x_1\tensor x_2)$ be a pointed plot of the balanced tensor product. Since $\pi_H$ is already a subduction, we can find a plot $\beta:V\to X\times_{H_0}\overline{X}$, defined on an open neighbourhood $0\in V\subseteq U_\alpha$ of the origin, such that $\alpha|_V=\pi_H\circ \beta$. This plot decomposes into two plots $\beta_1,\beta_2\in\diffeology_X$ on $X$, satisfying $r_X\circ\beta_1=r_X\circ\beta_2$. We use the notation $\alpha|_V=\beta_1\tensor \beta_2$. In particular, we get an equality $x_1\tensor x_2=\beta_1(0)\tensor \beta_2(0)$ inside the balanced tensor product, which means that we can find an arrow $h\in H$ such that $\beta_i(0)=x_ih$. The target must be $\target(h)=r_X(x_1)=r_X(x_2)$. This arrow allows us to write a pointed plot $r_X\circ\beta_i:(V,0)\to (H_0,\target(h^{-1}))$, so that now we can use that $\grpd{H}$ is locally subductive. Since the target map of $H$ is a local subduction, we can find a pointed plot $\Omega:(W,0)\to (H,h^{-1})$ such that $r_X\circ\beta_i|_W = \target_H\circ\Omega$. This relation means that, for every $t\in W$, we have a well-defined action $\beta_i(t)\cdot\Omega(t)\in X$. Hence we get a pointed plot
			\begin{equation*}
				\Psi:(W,0)\longrightarrow (X\times_{H_0}^{r_X,r_X}\overline{X},(x_1,x_2));
				\qquad
				t\longmapsto\left(\beta_1(t)\Omega(t),\beta_2(t)\Omega(t)\right).
			\end{equation*}
			It then follows by the definition of the balanced tensor product that
			\begin{equation*}
				\pi_H\circ \Psi(t) = \beta_1|_W(t)\Omega(t)\tensor \beta_2|_W(t)\Omega(t) = \beta_1|_W(t)\tensor \beta_2|_W(t) =\alpha|_W(t),
			\end{equation*}
			proving that $\pi_H$ is a local subduction.
		\end{proof}
	\end{lemma}
	
	\begin{lemma}
		\label{lemma:biprincipal bib between loc subd groupoids is locally bisubductive}
		If $G\LAalong{l_X}X\RAalong{r_X}H$ is a biprincipal bibundle between locally subductive groupoids, then the moment maps $l_X$ and $r_X$ are local subductions as well.
		
		\begin{proof}
			If the bibundle $G\LAalong{l_X}X\RAalong{r_X}H$ is biprincipal, we get two biequivariant diffeomorphisms $\varphi_G:X\tensor_H\overline{X}\to G$ and $\varphi_H:\overline{X}\tensor_G X\to H$ (\cref{proposition:biprincipal bibundle is weakly invertible}). It follows that the local subductivity of the source and target maps of $G$ and $H$ transfer to the four moment maps of the balanced tensor products. For example, the left moment map $L_X:X\tensor_H\overline{X}\to G_0$ can be written as $L_X= \target_G\circ\varphi_G$, where the right hand side is clearly a local subduction. We know as well that $L_X$ fits into a commutative square with the original moment map $l_X$:
			\begin{equation*}
				\begin{tikzcd}
					X\times_{H_0}^{r_X,r_X}\overline{X}\arrow[r,"\pi_H"]\arrow[d,"\pr_1|_{X\times_{H_0}\overline{X}}"'] & X\tensor_H \overline{X}\arrow[d,"L_X"]\\
					X\arrow[r,"l_X"'] & G_0.
				\end{tikzcd}
			\end{equation*}
			Since local subductions compose, and since by \cref{lemma:projection onto balanced tensor product is local subduction} the projection $\pi_H$ is a local subduction, we find that the upper right corner $L_X\circ\pi_H$ must be a local subduction. Hence the composition $l_X\circ \pr_1|_{X\times_{H_0}\overline{X}}$ is a local subduction, which by an argument that is analogous to the proof of \cref{lemma:properties of subductions}\emph{(2)} gives the local subductiveness of $l_X$. That the right moment map $r_X$ is a local subduction follows from a similar argument.
		\end{proof}
	\end{lemma}
	
	The lemma suggests that, if we refine our notion of principality something we might call \demph{pure-principality}, by passing from subductions to local subductions, then biprincipality between locally subductive groupoids means the same thing as this new notion of pure-principality. Let us make this precise.
	
	\begin{definition}
		Two diffeological groupoids are called \emph{purely Morita equivalent} if there exists a biprincipal bibundle between them, such that the two underlying moment maps are local subductions.
	\end{definition}
	
	Clearly, pure Morita equivalence implies ordinary Morita equivalence in the sense of \cref{definition:Morita equivalence and biprincipality}, since local subductions are, in particular, subductions. The question is if the converse implication holds as well. We have a partial answer, since \cref{lemma:biprincipal bib between loc subd groupoids is locally bisubductive} can now be restated as follows:
	\begin{proposition}
		\label{proposition:pure-ME is the same as ME for locally subductive groupoids}
		Two locally subductive groupoids are Morita equivalent if and only if they are purely Morita equivalent.
	\end{proposition}
	
	Especially in light of the existence of subductions that are not local subductions (see e.g. \cite[Exercise 61, p.60]{iglesias2013diffeology}), and the fact that the proof of \cref{lemma:biprincipal bib between loc subd groupoids is locally bisubductive} relies so heavily on the assumption that the groupoids are locally subductive, it seems that the ordinary diffeological Morita equivalence of \cref{definition:Morita equivalence and biprincipality} is not equivalent to pure-Morita equivalence in general. We do not, however, know of an explicit counter-example. This discussion leaves us an open question:
	\begin{question}
		\label{question:does inclusion pseudofunctor reflect weak equivalence}
		Does diffeological Morita equivalence reduce to Lie Morita equivalence on Lie groupoids? That is to ask, if two Lie groupoids are diffeologically Morita equivalent, are they also Lie Morita equivalent?
	\end{question}
	
	If two Lie groupoids $G$ and $H$ are diffeologically Morita equivalent, then there exists a diffeological biprincipal bibundle $G\LAalong{l_X}X\RAalong{r_X}H$, where $X$ is a diffeological space. A positive answer to \cref{question:does inclusion pseudofunctor reflect weak equivalence} could consist of a proof that $X$ is in fact a smooth manifold. Since $G$ and $H$ are both manifolds, it follows that $X\tensor_H\overline{X}$ and $\overline{X}\tensor_G X$ are also manifolds. We do not know if this is sufficient to imply that $X$ itself has to be a manifold. One suggestion is to use \cite[Article 4.6]{iglesias2013diffeology}, which gives a characterisation for when a quotient of a diffeological space by an equivalence relation is a smooth manifold. Since the balanced tensor products are quotients of diffeological spaces, one may try to use this result to obtain a special family of plots for their underlying fibred products. This could potentially be used to define an atlas on $X$.

	\subsection{Directions for future research}
	We list here some possible directions for future research. These are also proposed at the end of \cite[Section 1.2.3]{schaaf2020diffeology-groupoids-and-ME}.
	\begin{itemize}
		\item Finding an answer to the open \cref{question:does inclusion pseudofunctor reflect weak equivalence} about \emph{diffeological} Morita equivalence between \emph{Lie} groupoids.
		
		\item The construction of a theory of bibundles for a more general framework of generalised smooth spaces. One possibility is to look at the \emph{generalised spaces} of \cite[Definition 4.11]{baez2011convenient} (subsuming diffeology), or even to look at arbitrary classes of sheaves. What is the relation between our theory of Morita equivalence and the discussion in \cite{meyer2015groupoids}? A theory of principal bibundles seems to exist in a general setting for groupoids in $\infty$-toposes: \cite{nlab2018bibundle}.

		\item What is the precise relation between differentiable stacks and diffeological groupoids (cf. \cite{watts2019diffeological})? Using our notion of Morita equivalence, what types of objects are \emph{``diffeological stacks''} (i.e., Morita equivalence classes of diffeological groupoids)?
		
		\item Can the \emph{Hausdorff Morita equivalence} for holonomy groupoids of singular foliations introduced in \cite{garmendia2019hausdorff} be understood as a Morita equivalence between diffeological groupoids?
		
		\item Can the bridge between diffeology and noncommutative geometry that is being built in \cite{bertozzini2016spectral,iglesias2018noncommutative,androulidakis2019diffeological,iglesias2020quasifolds} be strengthened by our theory of Morita equivalence? Morita equivalence of Lie groupoids is already an important concept in relation to noncommutative geometry, especially for the theory of groupoid \Cstar-algebras. Can this link be extended to the diffeological setting, possibly through a theory of groupoid \Cstar-algebras for (a large class of) diffeological groupoids? If such a theory exists, what is the relation between Morita equivalence of diffeological groupoids and the Morita equivalence of their groupoid \Cstar-algebras? Is Morita equivalence preserved just like in the Lie case?
	\end{itemize}
	
	\printbibliography
	
\end{document}